\documentclass[11pt,letterpaper]{article}
\pdftrailerid{}
\usepackage[utf8]{inputenc}
\usepackage[letterpaper,margin=1in]{geometry}
\usepackage{amsmath,amssymb,amsthm,mathtools}
\usepackage{microtype,graphicx,enumitem,needspace,booktabs,tabularx}
\usepackage[skip=5pt,indent=0pt]{parskip}
\usepackage[section]{placeins}
\usepackage[numbers,sort&compress]{natbib}
\usepackage{xcolor}
\usepackage[unicode,bookmarksnumbered,bookmarksdepth=3,
  hidelinks]{hyperref}
\numberwithin{equation}{section}
\newtheorem{theorem}{Theorem}
\newtheorem{lemma}[theorem]{Lemma}
\newtheorem{proposition}[theorem]{Proposition}
\newtheorem{corollary}[theorem]{Corollary}
\setlist{itemsep=3pt,topsep=4pt}
\allowdisplaybreaks[1]
\newcommand{\paperauthors}{Zhou Longfei\thanks{\href{mailto:longfei.gabriel.zhou@gmail.com}{\texttt{longfei.gabriel.zhou@gmail.com}}}}

\title{Sharp Volume Contraction of a Centrally Cut John Ellipsoid}
\author{\paperauthors}
\date{}
\hypersetup{
  pdftitle={Sharp Volume Contraction of a Centrally Cut John Ellipsoid},
  pdfauthor={Zhou Longfei},
  pdfsubject={Sharp central-cut contraction, equality, and quantitative stability},
  pdfkeywords={John ellipsoid, convex geometry, convex optimization, cutting planes}}
\begin{document}
\maketitle
\begin{abstract}
For a convex body $K\subset\mathbb R^n$ in positive dimension, let $w(K)$ be
the volume of its maximum-volume inscribed ellipsoid. For a closed halfspace
$H$ whose boundary passes through the center of that ellipsoid, we prove
\[
w(K\cap H)\le R_nw(K),\qquad R_1=\frac12,\qquad
R_n=\frac12\left(\frac n{n-1}\right)^{(n-1)/2}\quad(n\ge2),
\]
with equality attained in every fixed dimension. This refines the optimal
dimension-independent bound from the author's earlier joint work on
Khachiyan's ellipsoid conjecture. We characterize the full equality family,
including nonsimplex extremizers, and establish an exact volume-deficit
identity that yields dimension-independent stability of the new ellipsoid
and its weighted contact geometry in John coordinates. The proof reduces
the problem to finite John contacts and constructs an explicit conic
certificate by a transformation of contact directions and a common polar
rotation. Its determinant is controlled by a nonnegative contact defect.
The sharp contraction theorem, equality classification, stability estimates,
and supporting results have been formally verified in Lean 4 using mathlib.

\end{abstract}
{\small
\par}
\section{Introduction}\label{sec:1}

A convex body is a compact convex subset of $\mathbb R^n$ with nonempty interior. Let $J(K)$ be its unique maximum-volume inscribed ellipsoid, let $c(K)$ be its center, and set
\begin{equation}\label{eq:1.1}
w(K)=\operatorname{vol}J(K)
=\max\{\operatorname{vol}E:E\subseteq K,\ E\text{ is an ellipsoid}\},
\end{equation}
where $\operatorname{vol}$ denotes Lebesgue volume. Existence, uniqueness, and the contact characterization are classical \cite{John1948,GG2007}. We study the change in $w$ under a halfspace cut whose boundary contains $c(K)$.

The volume of the recomputed ellipsoid provides a geometric measure of
progress after a cut. Determining its sharp decrease also raises a structural
question: which bodies attain the bound, and what does a small deficit reveal
about their ellipsoids and supporting contacts? We address these questions
through one explicit certificate and its determinant estimate.

\Needspace{8\baselineskip}
\subsection{Main results}\label{sec:1.2}

The results follow two connected threads.

\paragraph{Sharp contraction.} Let $C_n$ be the supremum of $w(K\cap H)/w(K)$ over convex bodies and halfspaces through their John centers. We prove
\begin{equation}\label{eq:1.3}
C_n=R_n,\qquad R_1=\frac12,\qquad
R_n=\frac12\left(\frac n{n-1}\right)^{(n-1)/2}\quad(n\ge2),
\end{equation}
with attainment in every dimension. The same sharp constant is obtained by restricting the supremum to finite John contact polytopes: the reduction in Section~\ref{sec:2.4} preserves the original ellipsoid volume while enlarging the cut region. The proof constructs an explicit conic certificate from the resulting finite data and bounds its determinant. A simple extremal configuration determines the candidate ellipsoid and the certificate parameters.

\begin{theorem}[{sharp central-cut contraction}]\label{thm-sharp}
Let $K\subset\mathbb R^n$, $n\ge1$, be a convex body, and let $H$ be a closed halfspace whose boundary passes through $c(K)$. Then
\begin{equation}\label{eq:3.1}
w(K\cap H)\le R_nw(K),
\end{equation}
where $R_n$ is given by \eqref{eq:1.3}. For every fixed dimension, there are bodies and central cuts attaining equality.

\end{theorem}

\paragraph{Equality and stability.} The logarithmic shortfall from $R_n$ is the sum of a geometric determinant deficit, a cone-complementarity term, and a positive-definite matrix divergence. Their zero sets give all equality cases, including nonsimplex extremizers. Their quantitative control yields dimension-independent bounds for the new ellipsoid's shape and center and for the weighted contact geometry.

For the following equality and stability statements, take John coordinates
$J(K)=B_2^n$, assume $n\ge2$, and write
$H_v=\{x:v^Tx\le0\}$ with $\|v\|=1$. Here $B_2^n$ is the Euclidean unit ball.
Put $b=1/n$ and $P=I-vv^T$. The extremal ellipsoid and contact set are
\[
\begin{gathered}
c_*=-\frac12v,\qquad
A_*=\frac12vv^T+\frac1{\sqrt{1-b}}P,\qquad
\mathcal E_*=c_*+A_*B_2^n,\\
\mathcal Z_b=\{-v\}\cup\{u:\|u\|=1,\ u^Tv=b\}.
\end{gathered}
\]
Theorem~\ref{thm-equality} states that equality holds if and only if
$\mathcal E_*\subseteq K$ and $K$ admits a supporting John frame whose
contacts lie in $\mathcal Z_b$. John frames are
the positive-weight contact decompositions defined in Section~\ref{sec:john-frames}.
At equality, $J(K\cap H_v)=\mathcal E_*$, and every positive-weight John frame
supporting $K$ has this pole-and-latitude geometry. The containment condition
leaves freedom in the rest of the body, giving extremizers beyond regular simplices.

For any feasible cut ellipsoid
$\mathcal E(c,A)=c+AB_2^n\subseteq K\cap H_v$, with $A$ symmetric positive
definite, suppose that its deficit satisfies
\[
0\le\epsilon:=\log\frac{R_n}{\det A}\le10^{-6}.
\]
Theorem~\ref{thm-stability} gives the explicit, dimension-independent estimates
\[
\|A-A_*\|_F\le30\sqrt\epsilon,\qquad
\|c-c_*\|\le50\sqrt\epsilon,\qquad
\sum_i\lambda_i\operatorname{dist}(u_i,\mathcal Z_b)^2
\le\frac{1024}{3}\epsilon,
\]
where the matrix norm is Frobenius and $(u_i,\lambda_i)$ is any supporting
John frame of $K$. Thus near equality controls the ellipsoid and the original
weighted contacts without selecting a unique extremal body. Section~\ref{sec:eq-stability}
proves the full classification and all stability estimates.

\paragraph{Certificates and optimization.}
On a bounded full-dimensional polytope $\mathcal P$ in any positive dimension,
the ellipsoid problem maximizes $\log\det A$ under second-order-cone support constraints.
For a normalized primal--dual feasible pair, with ellipsoid $\mathcal E(c,A)$
and positive-definite dual aggregate $M$, the computable gap satisfies
\[
G=-\log\det M-\log\det A\ge0,\qquad
\operatorname{vol}\mathcal E(c,A)\ge e^{-G}w(\mathcal P).
\]
The normalization and feasibility conditions are given in
Section~\ref{sec:2.2}; the relative-volume guarantee is established in
Section~\ref{sec:7.3}. It supplies a stopping criterion for an ellipsoid
subproblem. For central cuts in dimensions at least two, our explicit certificate bounds the new optimum
from the original John contacts, without requiring the unknown recomputed
ellipsoid. The exact deficit identity in Proposition~\ref{prop-gap} then
separates the geometric determinant loss from the optimization gap.
This connects the sharp contraction analysis with certifiable ellipsoid
computation. Exact duality and contact recovery are proved in
Section~\ref{sec:2.3}.

John's theorem and the finite convex optimization model are background. The problem-specific contribution is the quantitative certificate estimate and the geometric structure that it reveals.

\paragraph{Organization.}\label{sec:1.3}
Theorem~\ref{thm-sharp} is proved in Section~\ref{sec:3}: finite contact
reduction leads to an explicit certificate, whose determinant is bounded in
all dimensions before the geometric argument is completed. Section~\ref{sec:eq-stability}
then derives the equality classification and stability from the exact deficit.
Section~\ref{sec:further} develops exact duality, iteration bounds and further
matrix consequences. The conclusion is in Section~\ref{sec:11}.

Table~\ref{tab:main-results} summarizes the principal results and their optimization interpretation.

\begin{table}[htbp]
\centering
\small
\begin{tabularx}{\linewidth}{@{}>{\raggedright\arraybackslash}p{0.20\linewidth}>{\raggedright\arraybackslash}X>{\raggedright\arraybackslash}p{0.22\linewidth}@{}}
\toprule
Result & Mathematical content & Location \\
\midrule
Sharp contraction & Exact attained factor in every positive dimension & Theorem~\ref{thm-sharp} \\
Equality & Pole-and-latitude contacts, model containment and the full nonunique extremal family & Theorem~\ref{thm-equality} \\
Stability & Explicit shape, center and weighted-contact bounds with dimension-independent constants & Theorem~\ref{thm-stability} \\
Certificates and gap & Feasible volume bounds and certified relative-volume accuracy for ellipsoid computation & Lemma~\ref{lem-certificate};\newline Section~\ref{sec:7.3} \\
\bottomrule
\end{tabularx}
\caption{The principal results and their optimization interpretation.}\label{tab:main-results}
\end{table}

\subsection{Related work}\label{sec:1.1}

Such cuts arise in localization methods for convex optimization. An oracle supplies a halfspace containing the target set, and the current region is replaced by its intersection with that halfspace. The choice of center controls the geometric progress. In the centroid method, Grünbaum's theorem bounds the surviving ordinary-volume fraction by $1-(n/(n+1))^n<1-e^{-1}$ \cite{Grunbaum1960}, but computing the centroid is expensive. The classical enclosing-ellipsoid method admits an $O(n^2)$ arithmetic update, with a dimension-dependent contraction factor \cite{GLS1988,Todd2023}. Analytic and volumetric centers offer other computational tradeoffs \cite{BV2004,AV1995,Vaidya1996,Anstreicher1997}.

The method of inscribed ellipsoids instead uses $c(K)$ and the potential $w(K)$ \cite{TKE1988}. It also provides a geometric normalization. If $J(K)=c+AB_2^n$, John's rounding theorem gives
\begin{equation}\label{eq:1.2}
c+AB_2^n\subseteq K\subseteq c+nAB_2^n,
\qquad
B_2^n\subseteq A^{-1}(K-c)\subseteq nB_2^n.
\end{equation}

For a centrally symmetric body, the factor $n$ improves to $\sqrt n$. This rounding is useful in geometric preconditioning, volume approximation, and fixed-dimensional integer programming \cite{GG2007,KT1993}.

For a polytope with $N$ inequalities, computing $J(K)$ is a finite convex optimization problem; Section~\ref{sec:2} gives its formulation and dual certificates. Interior-point methods exploit its structure \cite{NN1994,KT1993,ZG2003}. One established bound computes an ellipsoid of volume at least $e^{-\delta}w(K)$, $0<\delta<1$, under a known normalization $B_2^n\subseteq K\subseteq\Lambda B_2^n$, within $O(N^{7/2}\log(N\Lambda/\delta))$ arithmetic operations \cite{Anstreicher2002}. This is a bound for approximate ellipsoid computation. Theorem~\ref{thm-sharp} concerns the exact John center; relating a computed center to an admissible cut requires the corresponding approximation analysis. Section~\ref{sec:7} distinguishes this computational cost from the number of exact central cuts.

Khachiyan asked for the optimal dimension-independent contraction factor \cite{Khachiyan1990}. Zhou, Zou, and Liu proved the value $\sqrt e/2$ and its uniform optimality \cite{Zhou2026}. The present paper determines the exact value in each dimension and describes the geometry of equality and near equality.

\section{Preliminaries}\label{sec:2}

We use affine normalization, John's contact characterization and a finite
volume certificate. These ingredients suffice for the sharp geometric bound.
The exact optimization interpretation of the certificate is developed later
in Section~\ref{sec:2.3}.

\subsection{Ellipsoids and finite optimization}\label{sec:2.1}

Throughout, the dimension is a positive integer. Let $B_2^n=\{z:\|z\|_2\le1\}$ be the closed Euclidean unit ball. Write
\begin{equation}\label{eq:2.1}
\mathcal E(c,A)=c+AB_2^n,\qquad A\in\mathbb S_{++}^n,
\qquad \operatorname{vol}\mathcal E(c,A)=\kappa_n\det A,
\quad \kappa_n=\operatorname{vol}(B_2^n).
\end{equation}

The matrix $A$ is the linear map applied to the ball. Its eigenvalues are the semiaxis lengths. For a single inequality $g^Tx\le\beta$, containment is equivalent to
\begin{equation}\label{eq:2.2}
\max_{\|z\|\le1}g^T(c+Az)=g^Tc+\|Ag\|\le\beta.
\end{equation}

Thus, for a bounded full-dimensional polytope
$\mathcal P=\{x:g_r^Tx\le\beta_r,\ r=1,\ldots,N\}$, the ellipsoid problem is
\begin{equation}\label{eq:2.3}
\begin{aligned}
\text{maximize}_{c,\ A\succ0}\quad&\log\det A\\
\text{subject to}\quad&g_r^Tc+\|Ag_r\|\le\beta_r,
\qquad r=1,\ldots,N.
\end{aligned}
\end{equation}

The objective is concave and the constraints are convex \cite[Section 8.4.2]{BV2004}. Each support constraint is also a second-order-cone constraint:
$(\beta_r-g_r^Tc,Ag_r)\in\mathcal L_{n+1}$, where
$\mathcal L_{n+1}=\{(\tau,z):\tau\ge\|z\|\}$.
We remove redundant inequalities with $g_r=0$.

\subsection{John frames and contact reduction}\label{sec:john-frames}

The following standard characterization turns maximality of the unit ball
into finite moment identities and supporting inequalities.

\begin{lemma}[{John's contact characterization \cite{John1948} and \cite[Theorems 3.1 and 3.3]{GG2007}}]\label{lem-john}
Let $K\subset\mathbb R^n$ be a convex body containing $B_2^n$.
Then $J(K)=B_2^n$ if and only if there are finitely many unit vectors
$u_i\in\partial K\cap S^{n-1}$ and positive weights $\lambda_i$ such that
\begin{equation}\label{eq:2.9}
u_i^Tx\le1\ (x\in K),\qquad
\sum_i\lambda_i u_i=0,\qquad
\sum_i\lambda_i u_i u_i^T=I,
\qquad \sum_i\lambda_i=n.
\end{equation}

\end{lemma}

We use this standard result without proof. The last identity follows by taking
traces in the second-moment identity. We call finite unit vectors and positive
weights with the two moment identities a \emph{John frame}; the inequalities in
\eqref{eq:2.9} say that the frame supports $K$. The same vector $u_i$ specifies both
the point of tangency with the unit ball and the outward normal of its
supporting plane.

\subsubsection{Finite contact reduction}\label{sec:2.4}

An invertible affine change of coordinates preserves ratios of ellipsoid
volumes. Send $J(K)$ to the unit ball and apply Lemma~\ref{lem-john}.
The resulting supporting frame converts the geometric problem into a finite
polytope problem while preserving the sharp supremum.

The frame does not need to describe the whole boundary of $K$. Retaining only
these supporting planes produces a polytope containing $K$, while John's
characterization certifies that its original maximal ellipsoid is unchanged.

\begin{lemma}[{finite contact reduction}]\label{lem-reduction}
Let $K$ be a convex body with
$J(K)=B_2^n$, and let $\mathcal U$ be a supporting John frame from Lemma~\ref{lem-john}.
Then the contact polytope defined below is a convex body satisfying
\begin{equation}\label{eq:2.18}
P_{\mathcal U}=\{x:u_i^Tx\le1\text{ for all }i\},\qquad
K\subseteq P_{\mathcal U},\qquad J(P_{\mathcal U})=B_2^n.
\end{equation}

In particular, $w(P_{\mathcal U})=w(K)=\kappa_n$. For every central halfspace
$H_v=\{x:v^Tx\le0\}$ with $\|v\|=1$, each ellipsoid in $K\cap H_v$ is also
feasible in $P_{\mathcal U}\cap H_v$, so
$w(K\cap H_v)\le w(P_{\mathcal U}\cap H_v)$.
Both cut sets have nonempty interior because they contain the corresponding
half of the unit ball.

\end{lemma}

\begin{proof}
The main point is boundedness: the supporting halfspaces could
otherwise define an unbounded set. We obtain a radius bound directly from
the two John moments. Since $\|u_i\|=1$, the polyhedron contains $B_2^n$.
For $x\in P_{\mathcal U}$, put $r=\|x\|$ and $s_i=u_i^Tx$. Then
$-r\le s_i\le1$, $\sum_i\lambda_i s_i=0$, and $\sum_i\lambda_i s_i^2=r^2$.
The product $(1-s_i)(r+s_i)$ is nonnegative and quadratic in $s_i$, so its
weighted sum can be evaluated using exactly these moments:
\begin{equation}\label{eq:2.19}
\begin{aligned}
0\le\sum_i\lambda_i(1-s_i)(r+s_i)
&=r\sum_i\lambda_i+(1-r)\sum_i\lambda_i s_i-\sum_i\lambda_i s_i^2\\
&=nr-r^2.
\end{aligned}
\end{equation}

Thus $r(r-n)\le0$ and $r\ge0$ imply $r\le n$. The finite intersection is closed and convex, so $B_2^n\subseteq P_{\mathcal U}\subseteq nB_2^n$ makes it a convex body. John's characterization gives $J(P_{\mathcal U})=B_2^n$, and the supporting inequalities give $K\subseteq P_{\mathcal U}$. 
\end{proof}

\textbf{Why the reduction preserves the sharp constant.} The enlargement can increase
the post-cut ellipsoid volume, but it leaves the pre-cut volume unchanged.
Consequently,
\begin{equation}\label{eq:2.20}
\frac{w(K\cap H_v)}{w(K)}
\le\frac{w(P_{\mathcal U}\cap H_v)}{w(K)}
=\frac{w(P_{\mathcal U}\cap H_v)}{w(P_{\mathcal U})}.
\end{equation}

For a particular $K$ the inequality can be strict. For the extremal problem,
however, the two classes have the same supremum. Every normalized convex body
admits such an enlargement. Conversely, any John frame defines a convex body
$P_{\mathcal U}$ with John ellipsoid $B_2^n$: the boundedness argument above
uses only the frame identities, and Lemma~\ref{lem-john} identifies its ellipsoid.
Thus contact polytopes themselves belong to the original class, and
\begin{equation}\label{eq:2.21}
C_n
=\sup_{\substack{\mathcal U\ \text{a John frame}\\\|v\|=1}}
\frac{w(P_{\mathcal U}\cap H_v)}{w(P_{\mathcal U})}
=\sup_{\substack{\mathcal U\ \text{a John frame}\\\|v\|=1}}
\frac{w(P_{\mathcal U}\cap H_v)}{\kappa_n}.
\end{equation}

The first direction gives an upper bound for every convex body; the converse
shows that no sharpness is lost. For each fixed frame and cut,
$w(P_{\mathcal U}\cap H_v)/\kappa_n$ is the maximum feasible determinant
in the finite ellipsoid problem \eqref{eq:2.3}, with supporting
inequalities $(g_i,\beta_i)=(u_i,1)$ and the one added cut $(g_0,\beta_0)=(v,0)$.
The geometric problem has therefore become a universal inequality for finite
contact data. It is this finite problem that the next section certifies.

\subsection{A finite volume certificate}\label{sec:2.2}

To bound the volume of every feasible ellipsoid, assign a scalar $\alpha_r$ and
a vector $q_r$ to each supporting inequality. With
$\operatorname{sym}X=(X+X^T)/2$, the complete certificate conditions are
\begin{equation}\label{eq:2.4}
\begin{gathered}
\sum_r\alpha_rg_r=0,\qquad \sum_r\alpha_r\beta_r=n,\\
\|q_r\|\le\alpha_r\quad\text{for every }r,\\
M=\sum_r\operatorname{sym}(q_rg_r^T)\succ0.
\end{gathered}
\end{equation}

Balance cancels the unknown center, and the second equality fixes the scale.
The norm inequality is the second-order-cone constraint
$(\alpha_r,q_r)\in\mathcal L_{n+1}$; it also implies $\alpha_r\ge0$.
The final condition collects the directional information into a positive-definite
matrix. Thus \eqref{eq:2.4} includes every constraint used by the dual problem below.

\begin{lemma}[{volume certificate}]\label{lem-certificate}
If \eqref{eq:2.4} holds, every feasible ellipsoid in \eqref{eq:2.3} satisfies
\begin{equation}\label{eq:2.5}
\operatorname{tr}(AM)\le n,\qquad
\det A\le\frac1{\det M}.
\end{equation}

\end{lemma}

\begin{proof}
The two conclusions require two steps: first turn the support
constraints into a trace bound, then apply scalar AM--GM to obtain a determinant
bound. Aggregate the support inequalities; balance cancels the center:
\begin{equation}\label{eq:2.6}
\begin{aligned}
\operatorname{tr}(AM)
&=\sum_r q_r^TAg_r
\le\sum_r\alpha_r\|Ag_r\|\\
&\le\sum_r\alpha_r(\beta_r-g_r^Tc)
=n-c^T\sum_r\alpha_rg_r=n.
\end{aligned}
\end{equation}

To pass from trace to determinant, use the symmetric positive-definite matrix
$Y=M^{1/2}AM^{1/2}$. Its trace and determinant are the quantities we need.
For its eigenvalues $\rho_1,\ldots,\rho_n>0$, scalar AM--GM gives
\begin{equation}\label{eq:2.7}
\det A\det M=\det Y=\prod_{j=1}^n\rho_j
\le\left(\frac1n\sum_{j=1}^n\rho_j\right)^n
=\left(\frac{\operatorname{tr}(AM)}n\right)^n\le1.
\end{equation}

Dividing by $\det M>0$ proves \eqref{eq:2.5}. 
\end{proof}

The matrix $M$ is a certificate: checking its finite constraints is enough to bound the optimum. A larger determinant gives a stronger bound. Optimizing the certificate leads to the normalized dual problem
\begin{equation}\label{eq:2.8}
\text{maximize }\log\det M
\quad\text{over }(\alpha_r,q_r)\text{ satisfying \eqref{eq:2.4}}.
\end{equation}

\section{Proof of the sharp contraction bound}\label{sec:3}

\label{sec:3.1}
We now prove Theorem~\ref{thm-sharp}. The equality model identifies the
parameters of the certificate. After proving its determinant estimate, a
common rotation supplies a feasible dual matrix and completes the volume
bound. The one-dimensional case is treated at the end.

For $n\ge2$, take John coordinates so that $J(K)=B_2^n$, and write
$H=H_v=\{x:v^Tx\le0\}$ with $\|v\|=1$. A candidate $\mathcal E(c,A)$ then has relative volume $\det A$. It is enough to prove $\det A\le R_n$ for every candidate in $K\cap H_v$. The cut body has nonempty interior: it contains, for example, the ball with center $-v/4$ and radius $1/8$.

\subsection{The equality model}\label{sec:3.2}

The certificate arguments in this section and the equality and stability results in Section~\ref{sec:eq-stability} assume $n\ge2$.

Put $b=1/n$ and let $P=I-vv^T$ be the orthogonal projection onto $v^\perp$. Consider the candidate
\begin{equation}\label{eq:3.2}
c_*=-\frac12v,\qquad
A_*=\frac12vv^T+\frac1{\sqrt{1-b}}P,
\qquad \det A_*=R_n.
\end{equation}

Its center moves away from the cutting plane by $1/2$. Its axial semiaxis is $1/2$, so it touches the cut; its transverse semiaxes expand to $(1-b)^{-1/2}$. The product of these lengths is the proposed sharp factor.

Take a regular simplex circumscribed about the unit ball, with unit facet normals $u_1,\ldots,u_{n+1}$. Its standard contact identities are
\begin{equation}\label{eq:3.3}
u_i^Tu_j=-\frac1n\ (i\ne j),\qquad
\sum_i u_i=0,\qquad
\sum_i u_i u_i^T=\frac{n+1}{n}I.
\end{equation}

Thus $\lambda_i=n/(n+1)$ are John weights. Choose $v=-u_1$, so the contact heights are $-1$ and $b$. The three required support values are
\begin{equation}\label{eq:3.4}
\begin{aligned}
(-v)^Tc_*+\|A_*(-v)\|&=\frac12+\frac12=1,\\
u_i^Tc_*+\|A_*u_i\|
&=-\frac b2+\sqrt{\frac{b^2}{4}+\frac{1-b^2}{1-b}}=1
&& (u_i^Tv=b),\\
v^Tc_*+\|A_*v\|&=-\frac12+\frac12=0.
\end{aligned}
\end{equation}

Hence $\mathcal E(c_*,A_*)\subseteq\Delta\cap H_v$ and
$w(\Delta\cap H_v)/w(\Delta)\ge\det A_*=R_n$.

\begin{figure}[tbp]
\centering
\includegraphics[width=0.98\linewidth]{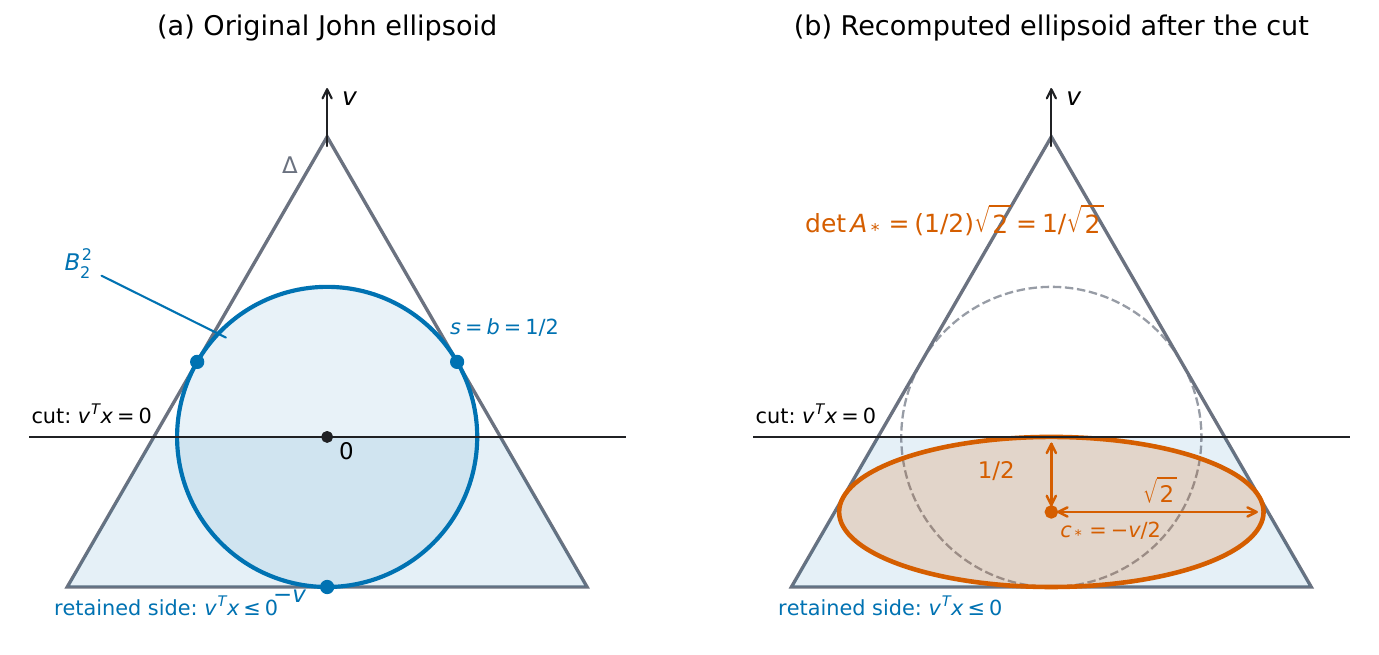}
\caption{The model in dimension two. The original inball has radius one. The new candidate has semiaxes $1/2$ and $\sqrt2$, so its area ratio is $1/\sqrt2$. The cut bisects the old ellipsoid, but the new ellipsoid can use space in the surrounding body.}\label{fig:1}
\end{figure}

The same support calculation applies to any John frame whose contacts have heights $-1$ and $b$. Section~\ref{sec:5} will describe all such equality configurations. For later use, write the inverse model shape as
$M_0=A_*^{-1}=2vv^T+\sqrt{1-b}\,P$.

\subsection{An explicit certificate from the contacts}\label{sec:3.3}

Let $(u_i,\lambda_i)$ be any supporting John frame. Define its axial and transverse components by $s_i=u_i^Tv$ and $y_i=Pu_i$. Old supporting planes carry the index $i$; the added cut has the separate label $0$.

The equality model tells us what the certificate should become:
$M_0=A_*^{-1}$. We seek formulas depending only on the original frame and the
cut, valid for every frame and matching that inverse shape at the model.
The scalar weights enforce balance; the vector directions enforce cone
feasibility and the model support relations.

\textbf{Balance the scalar weights.} Tilting the old weights linearly with the
contact height shifts their weighted normal sum in the direction $-v$.
The new cut, whose normal is $v$, can then restore balance. The John moments
evaluate both sums exactly. Choose
\begin{equation}\label{eq:3.5}
t=\frac{1+b}{2},\qquad
\alpha_i=\lambda_i(1-ts_i)\ge\lambda_i(1-t)\ge\frac{\lambda_i}{4},
\qquad \alpha_0=t.
\end{equation}

The John moments give normalization and balance directly:
\begin{equation}\label{eq:3.6}
\begin{aligned}
\sum_i\alpha_i
&=\sum_i\lambda_i-t\sum_i\lambda_i s_i=n,\\
\sum_i\alpha_i u_i+tv
&=\sum_i\lambda_i u_i-t\left(\sum_i\lambda_i u_iu_i^T\right)v+tv=0.
\end{aligned}
\end{equation}

The cut has right-hand side zero, so its weight does not enter the first sum.
The two identities are exactly the normalization and balance conditions in
\eqref{eq:2.4}. The bound $\alpha_i\ge\lambda_i/4$ keeps every original contact present
and will later control support slacks in the stability argument.
Balance alone allows other values of $t$; matching the equality model fixes
the value in \eqref{eq:3.5}, as explained below.

\textbf{Choose admissible vector directions.} We need unit directions for cone
feasibility: then the provisional choice $q_i=\alpha_i p_i$ lies on the
boundary of the cone. At equality these directions must also be
$A_*u_i/\|A_*u_i\|$, so that the support Cauchy--Schwarz inequality is exact.
The following transformation changes the axial coordinate while preserving
the unit sphere:
\begin{equation}\label{eq:3.7}
a=\frac b{2-b},\qquad k=\sqrt{1-a^2},\qquad
p_i=\frac{(s_i-a)v+ky_i}{1-as_i}.
\end{equation}

Since $\|y_i\|^2=1-s_i^2$ and $1-as_i\ge1-a\ge2/3$, the directions are well defined and
\begin{equation}\label{eq:3.8}
\|p_i\|^2
=\frac{(s_i-a)^2+(1-a^2)(1-s_i^2)}{(1-as_i)^2}=1.
\end{equation}

This is the unit-direction form of a Lorentz transformation. It fixes the
pole $-v$. At the model latitude, the required axial coordinate is
$b/(2+b)$, since $\|A_*u_i\|=1+b/2$ there. Matching
$(b-a)/(1-ab)=b/(2+b)$ gives the chosen $a=b/(2-b)$ and ensures
\begin{equation}\label{eq:3.9}
s_i\in\{-1,b\}\quad\Longrightarrow\quad
p_i=\frac{A_*u_i}{\|A_*u_i\|},
\end{equation}
the contact direction required at the model optimum. The transverse John
moment then shows that matching the transverse part of $M_0$ requires
$1-tb=(1-b)(1+b/2)$, which gives $t=(1+b)/2$.
Thus the parameters in \eqref{eq:3.5} and \eqref{eq:3.7} are fixed by exact model matching,
while the balance and unit-norm identities remain valid for arbitrary frames.

\textbf{The contact defect $\chi$.} To extend model matching to a quantitative
statement for arbitrary contacts, introduce
\begin{equation}\label{eq:3.10}
\chi=\sum_i\lambda_i\frac{(1+s_i)(s_i-b)^2}{1-as_i}.
\end{equation}

The factors $1+s_i$ and $(s_i-b)^2$ make each contribution vanish at the
model pole and latitude. The positive denominator is the same one appearing
in the transformed directions. This choice makes $\chi$ compatible with
the moment calculations of Section~\ref{sec:8}, while preserving
$\chi=0\iff s_i\in\{-1,b\}$ for every contact.

The defect depends only on the original weighted contact geometry and the
cut direction. It will control the determinant excess, identify all equality
frames, and bound weighted distances from those frames. It is therefore the
common geometric quantity behind the sharp bound, equality, and stability.

\textbf{The aggregate matrix $Q$.} Collect the original normals and their
transformed directions into one matrix:
\begin{equation}\label{eq:3.11}
Q=\sum_i\alpha_i u_ip_i^T+t\,vv^T.
\end{equation}

Each term $\alpha_i u_ip_i^T$ records one original constraint and its proposed
dual direction; the term $tvv^T$ represents the new cut.
Thus $Q^T$ is the aggregate of the unrotated vector multipliers.
Like $\chi$, it contains no unknown center or shape of the recomputed ellipsoid.

The following theorem supplies invertibility and a determinant lower bound.
Section~\ref{sec:3.4} will then use one common rotation to turn $Q^T$ into a
positive-definite certificate with the same determinant. This is how the finite
contact construction bounds every feasible cut ellipsoid without solving the
new optimization problem. The central estimate links its volume certificate,
$\log(R_n\det Q)$, to the contact defect $\chi$.

\begin{theorem}[{John-frame certificate bound}]\label{thm-frame-bound}
The matrix $Q$ is invertible with positive determinant. For $d=2-b$,
\begin{equation}\label{eq:3.12}
\frac{\chi}{1500}\le\log(R_n\det Q)\le\frac{\chi}{2d}.
\end{equation}

In particular, $\det Q\ge R_n^{-1}$. Equality in this determinant bound holds if and only if every contact has height $-1$ or $b$; then $Q=M_0$.

\end{theorem}

Theorem~\ref{thm-determinant} in Section~\ref{sec:8} is the normalized
form of this estimate. We prove it first, then return to the common polar
rotation and the geometric bound in Section~\ref{sec:3.4}.

\subsection{Factorization and the determinant estimate}\label{sec:8}

The remaining work is a determinant estimate for finite contact data.
We first factor the certificate through positive symmetric matrices. Scalar
interpolation treats higher dimensions, while transverse block information
handles dimensions two and three. The three branches are assembled in
Section~\ref{sec:8.6}.

\textbf{From the certificate to a symmetric matrix.} Recall $b=1/n$,
$P=I-vv^T$, and $u_i=s_iv+y_i$, where $s_i=u_i^Tv$ and $y_i=Pu_i$.
Use the parameters of Section~\ref{sec:3}:
$d=2-b$, $a=b/d$, $t=(1+b)/2$, $k=2\sqrt{1-b}/d$, and
$f(s)=(1-ts)/(1-as)$. These parameters depend only on the dimension; the
contact data enter through the weighted second-moment matrix $B$ below.
The useful factorization of the existing certificate is
\begin{equation}\label{eq:8.1}
\begin{aligned}
B&=\sum_i\lambda_i f(s_i)u_iu_i^T,&
F&=B+\frac{t}{1-a}vv^T,\\
D&=kP+k^2vv^T,&Q&=FD.
\end{aligned}
\end{equation}

Here $F$ adds a rank-one axial term to $B$, while $D$ rescales the transverse
directions by $k$ and the axial direction by $k^2$. This separates the
frame-dependent symmetric matrix $F$ from a fixed scaling whose determinant
is known.

To check the factorization, use
$\sum_i\lambda_i f(s_i)u_i=aBv-tv$ in the expansion of \eqref{eq:3.11}.
It gives $Q=BD+t(1+a)vv^T=FD$. All $f(s_i)$ are positive, so isotropy
makes $B$ and $F$ positive definite. Also $D\succ0$. This already proves
that $Q$ is invertible and has positive determinant.

\textbf{Choose the equality model as the reference.} At the model contacts, the
certificate is $Q=M_0=A_*^{-1}$. Thus the corresponding symmetric factor is
$F_0=M_0D^{-1}$. Normalize $F$ by this reference matrix:
\begin{equation}\label{eq:8.2}
\begin{gathered}
F_0=\frac{d^2}{2(1-b)}vv^T+\frac d2P,
\qquad\widehat F=F_0^{-1/2}FF_0^{-1/2},\\
F_0D=M_0,\qquad\det\widehat F=R_n\det Q.
\end{gathered}
\end{equation}

This congruence preserves positive definiteness and sends $F_0$ to the identity.
Consequently, $\widehat F-I$ measures the matrix departure from the equality
model. The last identity uses $\det D=k^{n+1}$ and
$\det(F_0D)=\det M_0=R_n^{-1}$. Thus the normalized determinant theorem
below is exactly the quantitative assertion needed in Theorem~\ref{thm-frame-bound}.

To measure the same departure directly in the contact data, recall the scalar
defect from \eqref{eq:3.10}:
\begin{equation}\label{eq:8.3}
\chi=\sum_i\lambda_i\frac{(1+s_i)(s_i-b)^2}{1-as_i}\ge0.
\end{equation}

The factors $1+s_i$ and $(s_i-b)^2$ make each summand vanish at the model pole
and latitude, respectively. The denominator is positive, and the John weights
retain the contribution of each contact. The next theorem links this scalar
description of contact geometry to the normalized determinant.

\begin{theorem}[{quantitative John-frame determinant inequality}]\label{thm-determinant}
Every John frame in dimension $n\ge2$ satisfies
\begin{equation}\label{eq:8.4}
\frac{\chi}{1500}\le\log\det\widehat F\le\frac{\chi}{2d},
\qquad
\det F\ge\det F_0=\frac{d^{n+1}}{2^n(1-b)}.
\end{equation}

Equality in the determinant inequality holds if and only if $s_i\in\{-1,b\}$ for every contact, in which case $F=F_0$ and $\widehat F=I$.

\end{theorem}

\textbf{Proof developed in Sections~\ref{sec:8.1}--\ref{sec:8.6}.} The dimension-specific lower
bounds are established in Sections~\ref{sec:8.2}, \ref{sec:8.3}, and \ref{sec:8.5}; Section~\ref{sec:8.6} assembles
them and proves the upper bound and equality statement.

The upper estimate will follow from the trace of $\widehat F-I$. The main work
is the lower estimate: a positive trace alone does not give the required determinant lower bound,
so we must use the relations between the matrix entries imposed by the contacts.
The coefficient $1/1500$ is the explicit margin produced by these estimates.
The sharp geometric constant is $R_n$; the proof does not assert optimality
of this quantitative margin.

\subsubsection{Scalar moment preparation}\label{sec:8.1}

For the high-dimensional argument, the needed contact information is the
height $s_i$ and the scalar defect $\chi$. Reweight by
$\nu_i=\lambda_i/(1-as_i)$ to absorb the denominator in $f$, and write
$N_r=\sum_i\nu_i s_i^r$. The scalar John moments give
\begin{equation}\label{eq:8.5}
N_0-aN_1=n,\qquad N_1-aN_2=0,\qquad N_2-aN_3=1.
\end{equation}

The three relations express $N_0,N_1,N_3$ in terms of $N_2$.
Expanding the defect gives
\begin{equation}\label{eq:8.6}
\chi=N_3+(1-2b)N_2+(b^2-2b)N_1+b^2N_0.
\end{equation}

Substitution expresses the axial entry of $B$ as
\begin{equation}\label{eq:8.7}
\zeta:=v^TBv=N_2-tN_3
=\frac{3d}{4}-\frac{d^2}{4(1-b)(2+b)}\chi.
\end{equation}

This is the scalar input used below: departure from the model lowers
$\zeta$ by a known multiple of $\chi$. All denominators are positive since
$0<b\le1/2$ and $-1\le s_i\le1$. The transverse information will be introduced
when the low-dimensional argument needs it.

\Needspace{10\baselineskip}

\subsection{The higher-dimensional estimate}\label{sec:8.2}

We now prove the lower bound in Theorem~\ref{thm-determinant} for $n\ge4$, so $0<b\le1/4$.
Two lemmas isolate the quantities to compare: Lemma~\ref{lem-lifted-bound} gives a determinant
gain and leaves a scalar logarithmic difference; Lemma~\ref{lem-hermite} bounds that
difference. Lemma~\ref{lem-rational} then completes the constant comparison.

\subsubsection{The lifted determinant bound}\label{sec:8.2.1}

Set $\ell(s)=\log f(s)$ and $\widetilde u_i=(u_i,\sqrt b)$.
Define the lifted moment matrix and weighted logarithmic sum by
\begin{equation}\label{eq:8.8}
\mathcal A=\sum_i\lambda_i f(s_i)\widetilde u_i\widetilde u_i^T,
\qquad \mathcal L=\sum_i\lambda_i\log f(s_i).
\end{equation}

The corresponding value for the pole-and-latitude model is
\begin{equation}\label{eq:8.9}
\mathcal L_0:=\frac{\ell(-1)}{1+b}+\frac{\ell(b)}{b(1+b)}.
\end{equation}

The two coefficients are the total John weights at the model pole and latitude.
Define also the gain coefficient
\begin{equation}\label{eq:8.10}
\delta_0(b):=\frac{(1+b)(b^2+b+2)}
{2d(1-b)(2+b)(3+b)}.
\end{equation}

\begin{lemma}[{lifted determinant bound}]\label{lem-lifted-bound}
For every John frame in dimension
$n\ge4$,
\begin{equation}\label{eq:8.11}
\log\det\widehat F
\ge (1+b)(\mathcal L-\mathcal L_0)+\delta_0(b)\chi.
\end{equation}

\end{lemma}

The matrix problem is thus reduced to a lower bound on
$\mathcal L-\mathcal L_0$. The next subsection estimates precisely this term.

\begin{proof}
\textbf{Encode the moments in the lifted matrix.} The constant last
coordinate turns the zero first moment into zero off-diagonal blocks, and
$b=1/n$ makes the final diagonal entry equal to one.

John's first and second moments combine into $\sum_i\lambda_i\widetilde u_i\widetilde u_i^T=I_{n+1}$. We use the standard finite-frame determinant estimate
\begin{equation}\label{eq:8.12}
XX^T=I_r,\quad \rho_i>0
\quad\Longrightarrow\quad
\log\det\bigl(X\operatorname{diag}(\rho_i)X^T\bigr)
\ge\sum_i\|X_i\|^2\log\rho_i,
\end{equation}
where $X_i$ are the columns. The point of this estimate is to replace a matrix
determinant by a weighted scalar sum. In the standard Cauchy--Binet/Jensen
argument \cite{HJ1985,BV2004}, the coefficients $w_S=(\det X_S)^2$, for subsets
$S$ of $r$ column indices, form a probability distribution. The needed
determinant identities are
\begin{equation}\label{eq:8.13}
\sum_S w_S=1,\qquad
\det(X\operatorname{diag}(\rho_i)X^T)=\sum_Sw_S\prod_{i\in S}\rho_i,
\qquad \sum_{S\ni i}w_S=\|X_i\|^2.
\end{equation}

For the marginal identity, expand the same determinant when only the $i$th diagonal weight is $1+z$:
$\det(I_r+zX_iX_i^T)=1+z\|X_i\|^2$.
Comparing coefficients of $z$ suffices. Scalar Jensen now gives
\begin{equation}\label{eq:8.14}
\log\!\left(\sum_Sw_S\prod_{i\in S}\rho_i\right)
\ge\sum_Sw_S\sum_{i\in S}\log\rho_i
=\sum_i\|X_i\|^2\log\rho_i,
\end{equation}
which proves \eqref{eq:8.12}. For $X_i=\sqrt{\lambda_i}\widetilde u_i$, this becomes
\begin{equation}\label{eq:8.15}
\log\det\mathcal A\ge(1+b)\mathcal L.
\end{equation}

\textbf{Return to the original determinant.} Both $\mathcal A$ and $F$ contain
$B$, so their determinant ratio can cancel this common matrix factor.
Write $\kappa=t/(1-a)$ for the axial coefficient already present in $F$,
and $\varrho=v^TB^{-1}v$ for the inverse axial moment introduced by the
Schur complement. The same moment identities used in \eqref{eq:8.1} give
\begin{equation}\label{eq:8.16}
\mathcal A=
\begin{pmatrix}
B&\sqrt b(aBv-tv)\\
\sqrt b(aBv-tv)^T&1-abt+a^2b\zeta
\end{pmatrix},
\qquad
\det\mathcal A=\frac{\det B}{n}(n+at-t^2\varrho).
\end{equation}

The Schur complement is positive. The abbreviation $H_0=n+at$ records its
constant part. Our aim is to remove $\varrho$ in favor of the direct axial
moment $\zeta=v^TBv$, whose dependence on $\chi$ is known from \eqref{eq:8.7}.
The determinant lemma and Cauchy--Schwarz give
\begin{equation}\label{eq:8.17}
\det Q=k^{n+1}\det B(1+\kappa\varrho),\qquad
\varrho\zeta\ge1,\qquad H_0-t^2\varrho>0.
\end{equation}

Cross multiplication, with positive denominators, therefore yields
\begin{equation}\label{eq:8.18}
\frac{\det Q}{\det\mathcal A}
=nk^{n+1}\frac{1+\kappa\varrho}{H_0-t^2\varrho}
\ge nk^{n+1}\frac{\zeta+\kappa}{H_0\zeta-t^2}.
\end{equation}

Taking logarithms and applying \eqref{eq:8.15}, we obtain
\begin{equation}\label{eq:8.19}
\log\det Q\ge(n+1)\log k+(1+b)\mathcal L+g(\zeta),
\quad
g(\zeta)=\log\frac{n(\zeta+\kappa)}{(n+at)\zeta-t^2}.
\end{equation}

All quantities inside the logarithm are positive: \eqref{eq:8.16} gives
$n+at>t^2\varrho\ge t^2/\zeta$.

\textbf{Extract the positive gain from the axial correction.} As $\chi$ increases,
\eqref{eq:8.7} makes $\zeta$ decrease from its model value $\zeta_0=3d/4$.
The variables $\theta_1,\theta_2$ below measure the relative decreases in the
numerator and denominator inside $g$. The denominator decreases more, which
increases their ratio. The elementary inequality $\log x\ge1-1/x$ will quantify
this increase. Write
\begin{equation}\label{eq:8.20}
\Delta_\zeta=\zeta_0-\zeta
=\frac{d^2\chi}{4(1-b)(2+b)},\qquad
\theta_1=\frac{\Delta_\zeta}{\zeta_0+\kappa},\qquad
\theta_2=\frac{H_0\Delta_\zeta}{H_0\zeta_0-t^2}.
\end{equation}

The positive-denominator condition implies $0\le\theta_1\le\theta_2<1$. Hence
\begin{equation}\label{eq:8.21}
\begin{aligned}
g(\zeta)-g(\zeta_0)
&=\log\frac{1-\theta_1}{1-\theta_2}
\ge1-\frac{1-\theta_2}{1-\theta_1}\\
&=\frac{\theta_2-\theta_1}{1-\theta_1}
\ge\theta_2-\theta_1=\delta_0(b)\chi,
\end{aligned}
\end{equation}

The last equality follows by substituting the parameters into the definition
of $\delta_0(b)$.

\textbf{Subtract the model value.} Since $f(-1)=d(3+b)/4$ and $f(b)=d/2$,
the model sum defined in \eqref{eq:8.9} satisfies
\begin{equation}\label{eq:8.22}
(n+1)\log k+(1+b)\mathcal L_0+g(\zeta_0)=\log R_n^{-1}.
\end{equation}

Combine \eqref{eq:8.19}, \eqref{eq:8.21}, and \eqref{eq:8.22}, and use
$\log\det\widehat F=\log R_n+\log\det Q$ from \eqref{eq:8.2}.
This gives \eqref{eq:8.11}, completing the proof. 
\end{proof}

\Needspace{8\baselineskip}

\subsubsection{The interpolation loss}\label{sec:8.2.2}

Keep $\ell(s)=\log f(s)$ and $\mathcal L=\sum_i\lambda_i\ell(s_i)$.
Recall explicitly the model value
$\mathcal L_0=\ell(-1)/(1+b)+\ell(b)/[b(1+b)]$ from \eqref{eq:8.9}.
Define the scalar integral
\begin{equation}\label{eq:8.23}
I(b):=\int_a^t\frac{z^2}{(1-z^2)(1-bz)^2}\,dz.
\end{equation}

\begin{lemma}[{interpolation loss}]\label{lem-hermite}
With these definitions, every John frame
in dimension $n\ge4$ satisfies
\begin{equation}\label{eq:8.24}
\mathcal L-\mathcal L_0\ge-(1-a)I(b)\chi.
\end{equation}

\end{lemma}

This is the scalar estimate needed in Lemma~\ref{lem-lifted-bound}. Its proof also gives the
exact pointwise remainder \eqref{eq:8.27}.

\begin{proof}
\textbf{Choose a quadratic with the correct weighted sum.} The John identities determine
the weighted average of every quadratic polynomial in the heights. We therefore
replace $\ell(s)=\log f(s)$ by a quadratic that matches the equality heights.
Matching the value at $-1$ and both the value and slope at $b$ also gives the
same vanishing factors as the defect $\chi$. Let $H_2$ satisfy
$H_2(-1)=\ell(-1)$, $H_2(b)=\ell(b)$, and $H_2'(b)=\ell'(b)$. Explicitly,
\begin{equation}\label{eq:8.25}
H_2(s)=\ell(b)+\ell'(b)(s-b)
+\frac{\ell(-1)-\ell(b)+(1+b)\ell'(b)}{(1+b)^2}(s-b)^2.
\end{equation}

The scalar John moments agree with those of masses $1/(1+b)$ at $-1$ and $1/[b(1+b)]$ at $b$. Since $H_2$ is quadratic,
\begin{equation}\label{eq:8.26}
\sum_i\lambda_iH_2(s_i)=\mathcal L_0.
\end{equation}

Therefore $\mathcal L-\mathcal L_0$ is the weighted sum of the remainders
$\ell(s_i)-H_2(s_i)$. We will establish the exact identity, for $-1\le s\le1$,
\begin{equation}\label{eq:8.27}
\ell(s)-H_2(s)=-(1+s)(s-b)^2
\int_a^t\frac{z^2\,dz}{(1+z)(1-bz)^2(1-zs)}.
\end{equation}

\textbf{Derive and bound the remainder.} We first prove a rational interpolation identity, then integrate
it to recover the logarithm. This exposes the sign and the vanishing factors
of the remainder explicitly.

\textbf{Interpolate the rational integrand.} Put $r_z(s)=s/(1-zs)$ and let
$H_{2,z}$ be the quadratic formula \eqref{eq:8.25} with data $r_z(-1)$, $r_z(b)$
and slope $(1-zb)^{-2}$. Clearing denominators gives the exact rational identity
\begin{equation}\label{eq:8.28}
r_z(s)-H_{2,z}(s)
=\frac{z^2(1+s)(s-b)^2}{(1+z)(1-bz)^2(1-zs)}.
\end{equation}

\textbf{Recover the logarithmic interpolant.} The value data integrate to
$\ell(-1)$ and $\ell(b)$ because $\ell(s)=-\int_a^t r_z(s)\,dz$.
It remains to check the slope data; these satisfy the elementary identity
\begin{equation}\label{eq:8.29}
-\int_a^t\frac{dz}{(1-bz)^2}
=-\frac{t-a}{(1-ab)(1-tb)}
=\frac{a}{1-ab}-\frac{t}{1-tb}=\ell'(b).
\end{equation}

Integration of the three polynomial coefficients therefore gives \eqref{eq:8.27}, with
no differentiation under the integral. All integrands are continuous and all
denominators are positive on the compact interval.

\textbf{Express the loss in units of the contact defect.} To sum the remainder,
we need a coefficient independent of the individual height $s$. The following
comparison leaves the factor defining $\chi$ outside the integral:
\begin{equation}\label{eq:8.30}
\frac{1-a}{1-z}-\frac{1-as}{1-zs}
=\frac{(z-a)(1-s)}{(1-z)(1-zs)}\ge0
\qquad(a\le z\le t,\ -1\le s\le1).
\end{equation}

It follows pointwise that
$\ell(s)-H_2(s)\ge-(1-a)I(b)(1+s)(s-b)^2/(1-as)$.
Summing with weights $\lambda_i$ proves \eqref{eq:8.24}. 
\end{proof}

\textbf{Combine the gain and loss estimates.} Lemmas~\ref{lem-lifted-bound} and \ref{lem-hermite} now give the
two inputs to the higher-dimensional part of Theorem~\ref{thm-determinant}. Substituting the
interpolation bound into \eqref{eq:8.11} yields
\begin{equation}\label{eq:8.31}
\log\det\widehat F\ge\delta_0(b)[1-\mathcal R(b)]\chi,
\qquad
\mathcal R(b)=\frac{4(1-b)^2(2+b)(3+b)}{b^2+b+2}I(b).
\end{equation}

Here $\mathcal R(b)=(1+b)(1-a)I(b)/\delta_0(b)$ is the loss as a fraction of
the gain. The remaining task is the scalar comparison in Lemma~\ref{lem-rational}.

\subsubsection{Exact evaluation and the scalar certificate}\label{sec:8.2.3}

The matrix argument is now complete except for a scalar comparison:
$\mathcal R(b)$ measures the interpolation loss as a fraction of the gain
$\delta_0(b)\chi$. We need a uniform margin below one and a positive lower bound
for $\delta_0$. The convenient targets below give
$(1/12)(1-124/125)=1/1500$, the coefficient in Theorem~\ref{thm-determinant}.

\begin{lemma}[{uniform scalar certificate}]\label{lem-rational}
For $0<b\le1/4$,
$\mathcal R(b)<124/125$ and $\delta_0(b)\ge1/12$.

\end{lemma}

\begin{proof}
We first reduce to the endpoint $b=1/4$, evaluate the integral by
partial fractions, and compare its exact value with $124/125$.
The lower bound for $\delta_0$ is a separate polynomial comparison at the end.

\textbf{Step 1: reduce to the largest parameter.} In the original integral, both
endpoints and the integrand depend on $b$. The substitution
$w=(z-a)/(1-az)$ fixes the lower endpoint at zero and exposes three factors
that can each be compared as $b$ varies:
\begin{equation}\label{eq:8.32}
\mathcal R(b)=\frac{4(2+b)(3+b)}{b^2+b+2}
\int_0^{(2+b)/(4+b)}\frac1{1-w^2}
\left(\frac{2w+b(1-w)}{2+b(1-w)}\right)^2dw.
\end{equation}

Temporarily name the outside prefactor $\pi$, the rational factor inside the
integrand $\psi_b$, and the upper endpoint $\ell_b$:
$\pi(b)=4(2+b)(3+b)/(b^2+b+2)$,
$\psi_b(w)=[2w+b(1-w)]/[2+b(1-w)]$, and
$\ell_b=(2+b)/(4+b)$. The largest endpoint is $\ell_{1/4}=9/17$.
For $0\le b_1\le b_2\le1/4$ and $0\le w\le9/17$,
\begin{equation}\label{eq:8.33}
\begin{aligned}
\pi(b_2)-\pi(b_1)
&=\frac{16(b_2-b_1)(1-b_1-b_2-b_1b_2)}
{(b_1^2+b_1+2)(b_2^2+b_2+2)}\ge0,\\
\psi_{b_2}(w)-\psi_{b_1}(w)
&=\frac{2(b_2-b_1)(1-w)^2}
{[2+b_1(1-w)][2+b_2(1-w)]}\ge0,\\
\ell_{b_2}-\ell_{b_1}
&=\frac{2(b_2-b_1)}{(4+b_1)(4+b_2)}\ge0.
\end{aligned}
\end{equation}

Here $1-b_1-b_2-b_1b_2\ge7/16$ and $\psi_b(w)\ge0$. The integrand and interval in \eqref{eq:8.32} are therefore increasing in $b$, giving
$\mathcal R(b)\le\mathcal R(1/4)$.

\textbf{Step 2: evaluate the endpoint integral by partial fractions.}
Specializing the same $w$-integral \eqref{eq:8.32} to $b=1/4$ gives
\begin{equation}\label{eq:8.34}
\mathcal R(1/4)=\frac{468}{37}
\int_0^{9/17}\frac{(1+7w)^2}{(1-w^2)(9-w)^2}\,dw.
\end{equation}

Its integrand has the elementary decomposition
\begin{equation}\label{eq:8.35}
\begin{aligned}
\frac{(1+7w)^2}{(1-w^2)(9-w)^2}
&=\frac{9}{50(1+w)}+\frac{1}{2(1-w)}\\
&\quad-\frac{8}{25(9-w)}-\frac{256}{5(9-w)^2}.
\end{aligned}
\end{equation}

All logarithm arguments are positive on $0\le w\le9/17$.
Integrating the four terms and collecting the endpoint logarithms yields
\begin{equation}\label{eq:8.36}
\begin{aligned}
\mathcal R(1/4)
&=\frac{468}{37}
\left(\frac9{50}\log13-\frac1{25}\log2-\frac{16}{45}\right)\\
&=\frac{2106\log13-468\log2-4160}{925}<\frac{124}{125}.
\end{aligned}
\end{equation}

Together with Step 1, this proves $1-\mathcal R(b)>1/125$ throughout the
required interval.

\textbf{Step 3: bound the gain coefficient.} We now check the other input to \eqref{eq:8.31}.
Clearing the positive denominators in $\delta_0(b)\ge1/12$ gives
\begin{equation}\label{eq:8.37}
6(1+b)(b^2+b+2)-d(1-b)(2+b)(3+b)
=b(26+19b+4b^2-b^3)\ge0.
\end{equation}

This proves both assertions. 
\end{proof}

\textbf{Higher-dimensional part of Theorem~\ref{thm-determinant} completed.} Combining \eqref{eq:8.31} with
Lemma~\ref{lem-rational} gives $\log\det\widehat F\ge\chi/1500$ for $n\ge4$.
The remaining dimensions follow below; the common upper bound and equality
statement are completed in Section~\ref{sec:8.6}.

\textbf{Why dimensions two and three require another estimate.} The scalar formula for
$\mathcal R(b)$ is well defined for $0<b\le1/2$. The preceding bound needs
$\mathcal R(b)<1$; for $n=3$ and $n=2$, however, this scalar remainder exceeds one.
Indeed, on $a\le z\le t$,
\begin{equation}\label{eq:8.38}
\frac1{1-z^2}\ge1+z^2+z^4,\qquad
\frac1{(1-bz)^2}\ge1+2bz+3b^2z^2.
\end{equation}

These follow by multiplying the positive denominators: the two products are
$1-z^6\le1$ and $1-(bz)^3(4-3bz)\le1$. With
$p(b)=4(1-b)^2(2+b)(3+b)/(b^2+b+2)$, exact polynomial integration gives
\begin{equation}\label{eq:8.39}
\begin{aligned}
\mathcal R(b)&\ge p(b)\int_a^t
z^2(1+z^2+z^4)(1+2bz+3b^2z^2)\,dz,\\
\mathcal R(1/3)&>\frac{16}{15}>1,\\
\mathcal R(1/2)&>\frac76>1.
\end{aligned}
\end{equation}

Thus the scalar estimate does not prove a positive determinant deficit in these
two dimensions. We now retain transverse block information, first to settle
dimension two and then to derive the additional control needed in dimension three.

\subsection{Transverse blocks and dimension two}\label{sec:8.3}

The scalar estimate has reached its limitation. To handle the remaining
dimensions, we retain the components of $u_i=s_iv+y_i$ transverse to $v$.
Work in $\mathbb Rv\oplus v^\perp$, with $m=n-1$.
In the normalized certificate, $x$ measures the axial decrease, $E$ is the
transverse perturbation, and $j$ couples the two subspaces.

\begin{lemma}[{block identities and constraints}]\label{lem-moments}
The normalized matrix satisfies
\begin{equation}\label{eq:8.40}
\widehat F=\begin{pmatrix}1-x&j^T\\j&I_m+E\end{pmatrix},
\quad x=\frac{\chi}{2(2+b)},\quad
\sigma:=\operatorname{tr}E=\frac{2\chi}{d(2+b)},\quad x=\frac d4\sigma,
\end{equation}
where
\begin{equation}\label{eq:8.41}
E=-\frac1d\sum_i\nu_i(s_i-b)y_i y_i^T,
\qquad
j=-\frac{\sqrt{1-b}}{b\sqrt d}\sum_i\nu_i(s_i-b)y_i.
\end{equation}

The scalar and matrix constraints are
\begin{equation}\label{eq:8.42}
0\le\sigma\le\frac{2(1-b)}{d+b^2},\qquad 0\le\chi\le2,
\end{equation}
\begin{equation}\label{eq:8.43}
\begin{pmatrix}x&-j^T\\-j&tI_m-E\end{pmatrix}\succeq0,
\qquad -\frac12I_m\prec E\prec tI_m.
\end{equation}

\end{lemma}

The trace of $E$ is nonnegative, but $E$ need not be positive semidefinite.
The Gram constraint limits the coupling: for every transverse vector $q$,
$(q^Tj)^2\le x\,q^T(tI_m-E)q$. At the equality contacts, either $y_i=0$
or $s_i-b=0$, so $E=j=0$.

\begin{proof}
\textbf{Recover the blocks.} The projected John identities are
$\sum_i\lambda_i y_i=\sum_i\lambda_i s_i y_i=0$ and
$\sum_i\lambda_i y_iy_i^T=I_m$. With the reweighted vector and matrix moments
$C_r=\sum_i\nu_i s_i^r y_i$ and $Y_r=\sum_i\nu_i s_i^r y_iy_i^T$, they give
\begin{equation}\label{eq:8.44}
C_0=aC_1,\qquad C_1=aC_2,\qquad Y_0-aY_1=I_m.
\end{equation}

The identity $f(s)=d/2-(s-b)/[2(1-as)]$ gives the transverse expression
for $E$, or equivalently $Y_1-bY_0=-dE$.
Adding the axial term $t/(1-a)$ to \eqref{eq:8.7} and normalizing by $F_0$ gives
$1-x$. Taking the trace of $E$, using $\|y_i\|^2=1-s_i^2$, gives $\sigma$.
The mixed entry follows from $C_1-tC_2$ and \eqref{eq:8.44}.
This proves \eqref{eq:8.40}--\eqref{eq:8.41}.

\textbf{Bound the scalar range.} It suffices to bound $\chi$.
The nonnegative sum $\sum_i\nu_i(1-s_i)(s_i+b)^2$ expands to
$-N_3+(1-2b)N_2+(2b-b^2)N_1+b^2N_0$.
Combining it with \eqref{eq:8.5}--\eqref{eq:8.6} cancels the remaining free moment:
\begin{equation}\label{eq:8.45}
\chi+\frac{(1-b)(2+b)^2}{(d+b^2)^2}
\sum_i\nu_i(1-s_i)(s_i+b)^2
=\frac{d(2+b)(1-b)}{d+b^2}.
\end{equation}

Nonnegativity gives the stated bound on $\sigma$ through \eqref{eq:8.40}.
The estimate $\chi\le2$ also uses
$2(d+b^2)-d(2+b)(1-b)=b(2+3b-b^2)\ge0$.

\textbf{Keep the coupling in a positive block matrix.} The required matrix has
the Gram representation
\begin{equation}\label{eq:8.46}
\begin{pmatrix}x&-j^T\\-j&tI_m-E\end{pmatrix}
=\sum_i\nu_i(1+s_i)r_ir_i^T,\qquad
r_i=\binom{(s_i-b)/\sqrt{2(2+b)}}
{\sqrt{(1-b)(2+b)/(2d)}\,y_i}.
\end{equation}

Its weights are nonnegative, proving the first part of \eqref{eq:8.43}.
Solving the two relations for $Y_0,Y_1$ gives
\begin{equation}\label{eq:8.47}
\begin{aligned}
\sum_i\nu_i(1-s_i)y_i y_i^T&=\frac d{2+b}(I_m+2E),\\
\sum_i\nu_i(1+s_i)y_i y_i^T&=\frac{2d}{(1-b)(2+b)}(tI_m-E).
\end{aligned}
\end{equation}

Both sums are positive definite: the contacts at $s_i=\pm1$ have $y_i=0$,
and all other transverse vectors span $v^\perp$ by the John second moment.
Thus $I_m+2E\succ0$ and $tI_m-E\succ0$, proving the strict spectral bounds. 
\end{proof}

\textbf{Dimension two.} Here $m=1$, so the trace determines the transverse block:
$E=\sigma$, $x=3\sigma/8$, and $t=3/4$.
The remaining coupling loss is $j^2$. Substituting the Gram bound
$j^2\le x(t-\sigma)$ into the determinant gives
\begin{equation}\label{eq:8.48}
\det\widehat F\ge(1-x)(1+\sigma)-x(t-\sigma)
=1+\frac{11}{32}\sigma=1+\frac{11}{60}\chi.
\end{equation}

Dimension three needs further control because the transverse trace alone no
longer determines the transverse determinant. We next state the precise
estimate needed there.

\subsection{A direct coupling inequality}\label{sec:8.4}

For the two-dimensional transverse block in the next subsection, we need to
control its deformation and its coupling to the axis together. The same
bound will also be used for local stability in Section~\ref{sec:9.3}.
The matrix $K_E=I_m+2E$ is positive definite by Lemma~\ref{lem-moments}.
The nonnegative quantity $A_E=\operatorname{tr}(E^2K_E^{-1})$ measures
transverse deformation: each eigenvalue $e$ of $E$ contributes $e^2/(1+2e)$.

\begin{lemma}[{direct coupling bound}]\label{lem-variational}
Let
$K_E=I_m+2E$, $A_E=\operatorname{tr}(E^2K_E^{-1})$, and
$L=d(1-b)/(d+b^2)-d\sigma/2\ge0$.
If $\sigma>0$, then
\begin{equation}\label{eq:8.49}
A_E\le\frac\sigma2,\qquad
\|K_E^{-1}j\|^2\le L\left(\frac\sigma2-A_E\right).
\end{equation}

\end{lemma}

The first inequality bounds transverse deformation. The second says that
the remaining amount $\sigma/2-A_E$ controls the transformed coupling.
The coefficient $L$ is the distance to the upper trace bound in
\eqref{eq:8.42}, multiplied by $d/2$.

\begin{proof}
\textbf{Evaluate only the tests needed for these conclusions.} We will take $Z=EK_E^{-1}$, first with $q=0$ and then with
$q=K_E^{-1}j/L$ when $L>0$. For these calculations set
\begin{equation}\label{eq:8.50}
\ell=\operatorname{tr}(ZE)+q^Tj,\qquad
Q_2=\operatorname{tr}(Z^2K_E)+4q^TZj+L\|q\|^2.
\end{equation}

To bound this expression directly from the frame, put
$\gamma=\sqrt{d(1-b)}/(d+b^2)$. The mixed moments, transverse identity,
and scalar relation in Section~\ref{sec:8.3} give, respectively,
\begin{equation}\label{eq:8.51}
\begin{aligned}
\sum_i\nu_i(1-s_i)y_iy_i^T&=\frac d{2+b}K_E,\\
\gamma\sum_i\nu_i(1-s_i)(s_i+b)y_i&=\frac{2d}{2+b}j,\\
\gamma^2\sum_i\nu_i(1-s_i)(s_i+b)^2&=\frac d{2+b}L.
\end{aligned}
\end{equation}

Together with the expression for $E$, these identities evaluate the linear
form and its associated squared norm:
\begin{equation}\label{eq:8.52}
\begin{aligned}
\ell&=-\frac1d\sum_i\nu_i(s_i-b)y_i^T
\bigl(Zy_i+\gamma(s_i+b)q\bigr),\\
\sum_i\nu_i(1-s_i)\|Zy_i+\gamma(s_i+b)q\|^2
&=\frac d{2+b}Q_2.
\end{aligned}
\end{equation}

Weighted Cauchy--Schwarz, with
$\|y_i\|^2=(1-s_i)(1+s_i)$, therefore gives
$\ell^2\le\chi Q_2/[d(2+b)]=(\sigma/2)Q_2$.
Pole terms have $y_i=0$, so this calculation requires no division by
$1\pm s_i$. Completing one scalar square yields
\begin{equation}\label{eq:8.53}
2\ell-Q_2\le2\ell-\frac{2\ell^2}{\sigma}
=\frac\sigma2-\frac2\sigma\left(\ell-\frac\sigma2\right)^2
\le\frac\sigma2.
\end{equation}

\textbf{Extract the deformation and coupling bounds.}
The chosen $Z=EK_E^{-1}$ is symmetric and makes
$\operatorname{tr}(ZE)=\operatorname{tr}(Z^2K_E)=A_E$.
With $q=0$, this gives $2\ell-Q_2=A_E$, proving the first bound.
Keeping this $Z$, use $2EK_E^{-1}=I_m-K_E^{-1}$: the part depending on $q$
is $2q^TK_E^{-1}j-L\|q\|^2$. For $L>0$ it is maximized at the chosen
$q=K_E^{-1}j/L$, giving
\begin{equation}\label{eq:8.54}
2\ell-Q_2=A_E+\frac1L\|K_E^{-1}j\|^2\le\frac\sigma2.
\end{equation}

If $L=0$, the same finite Cauchy--Schwarz calculation with $Z=0$ and $q=j$
gives $\|j\|^4\le0$, hence $j=0$. This proves the second bound at that endpoint.
The case $\sigma=0$ has $E=j=0$ by \eqref{eq:8.40}--\eqref{eq:8.41}, so the estimates extend
there as well. 
\end{proof}

The useful output is \eqref{eq:8.49}. Section~\ref{sec:10.4} records the more general inequality
for arbitrary tests; the proof here has already supplied everything needed
for the three-dimensional and local stability arguments.

\subsection{The three-dimensional estimate}\label{sec:8.5}

The transverse block is now two-dimensional, so its trace no longer determines
its determinant. We must control both the transverse determinant and its
interaction with $j$. The proof first expresses that interaction through two
scalar quadratic forms, then combines the Gram and variational inequalities
to remove them. The remaining task is a polynomial bound on an interval.

\begin{proposition}[{three-dimensional bound}]\label{prop-dim3}
In dimension three,
\begin{equation}\label{eq:8.55}
\det\widehat F\ge1+\frac2{11}\sigma
=1+\frac{36}{385}\chi.
\end{equation}

\end{proposition}

\begin{proof}
The desired output is a positive determinant excess proportional
to $\sigma$. The coefficient $2/11$ below is the margin supplied by the
polynomial estimate.

\textbf{Step 1: describe the transverse block by scalar invariants.} In dimension
three, the definitions specialize to $E\in\mathbb S^2$, $x=5\sigma/12$,
$t=2/3$, and $0\le\sigma\le3/4$. If $\sigma=0$, the defect vanishes and
\eqref{eq:8.41} gives $\widehat F=I$. Suppose $\sigma>0$ and write
\begin{equation}\label{eq:8.56}
p=\det E,\qquad z=\sigma+4p,\qquad
L=\frac58-\frac56\sigma,\qquad
\Delta=\frac{11}{3}(1+\sigma)-z.
\end{equation}

Here $p$ is the second invariant of $E$ besides its trace. The combination
$z=\sigma+4p$ is chosen because it expresses the variational remainder
$\sigma-2A_E$ by a simple fraction. The coefficient $L$ is exactly the one
in Lemma~\ref{lem-variational} at $n=3$. The quantity $\Delta$ will be the positive denominator
when we combine the two coupling inequalities.
For a two-by-two matrix, Cayley--Hamilton gives
\begin{equation}\label{eq:8.57}
\det K_E=1+\sigma+z,\qquad
\sigma-2A_E=\frac{z}{1+\sigma+z}.
\end{equation}

Lemma~\ref{lem-variational} and $\det E\le\sigma^2/4$ therefore imply
\begin{equation}\label{eq:8.58}
0\le z\le\sigma(1+\sigma),\qquad
j^TK_E^{-2}j\le\frac{Lz}{2(1+\sigma+z)},\qquad \Delta>0.
\end{equation}

Thus $z$ lies in an explicit interval once $\sigma$ is fixed. This will later
make a concavity argument possible.

\textbf{Step 2: bound exactly the coupling term in the determinant.} For a
two-by-two transverse block, that term is
$j^T\operatorname{adj}(I_2+E)j=(1+\sigma)\|j\|^2-j^TEj$.
Use the short names $r=\|j\|^2$ and $\omega=j^TEj$ for these two quadratic
forms, and set $N=tI_2-E\succ0$. The Gram constraint and \eqref{eq:8.58} read, respectively,
\begin{equation}\label{eq:8.59}
(t-\sigma)r+\omega\le x\det N,
\end{equation}
\begin{equation}\label{eq:8.60}
\bigl[(1+\sigma)(1+4\sigma)-z\bigr]r-4(1+\sigma)\omega
\le\frac L2z(1+\sigma+z).
\end{equation}

Both inequalities are linear in $r,\omega$. Choose a positive combination
whose coefficients are exactly those of $(1+\sigma)r-\omega$.
Solving this two-variable coefficient-matching problem gives the multipliers
$[3(1+\sigma)+z]/\Delta$ and $5/(3\Delta)$, respectively.
Their positivity preserves the inequality direction, and their common
denominator explains the definition of $\Delta$. Adding yields
\begin{equation}\label{eq:8.61}
j^T\operatorname{adj}(I_2+E)j
\le\frac{x\det N[3(1+\sigma)+z]+\frac56Lz(1+\sigma+z)}\Delta.
\end{equation}

\textbf{Step 3: eliminate the matrix variables.} Substitute this bound into
$\det\widehat F=(1-x)(1+\sigma+p)-j^T\operatorname{adj}(I_2+E)j$
and use $p=(z-\sigma)/4$. Multiplication by $108\Delta$ clears the variable
denominator and gives a polynomial in $\sigma,z$ alone; $108$ is a convenient
positive normalization of its rational coefficients:
\begin{equation}\label{eq:8.62}
108\Delta(\det\widehat F-1)\ge\mathcal P_*(\sigma,z),
\end{equation}
\begin{equation}\label{eq:8.63}
\mathcal P_*(\sigma,z)=72\sigma(1+\sigma)
+\left(\frac{171}{4}-\frac{53}{4}\sigma+75\sigma^2\right)z
+\left(-\frac{333}{4}+75\sigma\right)z^2.
\end{equation}

\textbf{Step 4: use concavity to check only the endpoints.} The coefficient of
$z^2$ is negative for $0\le\sigma\le3/4$, so the minimum over the admissible
$z$ interval occurs at an endpoint. At $z=0$ the polynomial is
$72\sigma(1+\sigma)$. We show that the other endpoint is at least this large.
At $z=\sigma(1+\sigma)$, it equals
\begin{equation}\label{eq:8.64}
\frac{\sigma(1+\sigma)}4
\bigl(300\sigma^3+267\sigma^2-386\sigma+459\bigr)
\ge72\sigma(1+\sigma),
\end{equation}
because the cubic in parentheses minus $288$ has the nonnegative decomposition
\begin{equation}\label{eq:8.65}
300\sigma^3+267\left(\sigma-\frac34\right)^2
+\frac{29}{2}\sigma+\frac{333}{16}.
\end{equation}

Consequently $\mathcal P_*\ge72\sigma(1+\sigma)$ throughout the admissible
interval. Returning through \eqref{eq:8.62}, the bound
$\Delta\le11(1+\sigma)/3$ cancels $1+\sigma$ and gives
$\det\widehat F-1\ge2\sigma/11$.
Thus $2/11=72/[108(11/3)]$ is the final margin left by the polynomial bound. 
\end{proof}

Keeping $L=5/8-5\sigma/6$ is decisive: a larger coefficient weakens the second constraint and loses the single-interval polynomial certificate. The proof uses only the moment identities, positive-semidefinite blocks, two-by-two matrix identities, and polynomial inequalities with rational coefficients.

\subsection{Completion of the determinant and geometric bounds}\label{sec:8.6}

\label{sec:3.4}
The three branches now supply determinant excess. We put them on the same
logarithmic scale, then recover the upper bound and the equality condition.

\textbf{Obtain one lower bound for all dimensions.}
For $n\ge4$, Lemma~\ref{lem-rational} and \eqref{eq:8.31} give $\log\det\widehat F\ge\chi/1500$.
For $n=2,3$, equations \eqref{eq:8.48} and \eqref{eq:8.55} give $\det\widehat F\ge1+c\chi$ with $c=11/60$ and $c=36/385$, respectively. Since $\chi\le2$ and $\log(1+u)\ge u/(1+u)$ for $u\ge0$,
\begin{equation}\label{eq:8.66}
\log\det\widehat F\ge\frac{c}{1+2c}\chi\ge\frac{\chi}{1500}.
\end{equation}

\textbf{Use the trace for the upper bound.} This part is common to all dimensions:
\begin{equation}\label{eq:8.67}
\log\det\widehat F\le\operatorname{tr}\widehat F-n
=\sigma-x=\frac{\chi}{2d}.
\end{equation}

\textbf{Recover equality from zero defect.} The lower bound forces $\chi=0$ when
$\det\widehat F=1$. Every summand in \eqref{eq:8.3} is nonnegative with a positive weight
and denominator, so this is equivalent to $s_i\in\{-1,b\}$ for all $i$.
Conversely, this contact condition gives $x=j=E=0$ by \eqref{eq:8.40}--\eqref{eq:8.41}.
This completes the proof of Theorem~\ref{thm-determinant} announced at the start of Section~\ref{sec:8}.

By \eqref{eq:8.2}, this proves \eqref{eq:3.12}. At equality $F=F_0$, hence $Q=M_0$; together with the positivity established above, this completes Theorem~\ref{thm-frame-bound}.

\paragraph{The common polar rotation.}
The matrix $Q^T$ need not be symmetric. A single rotation of all the directions produces the positive-definite aggregate required by Lemma~\ref{lem-certificate}. By Theorem~\ref{thm-frame-bound}, the left polar decomposition is well defined:
\begin{equation}\label{eq:3.13}
M=(QQ^T)^{1/2},\qquad U=M^{-1}Q,
\qquad q_i=\alpha_iUp_i,\qquad q_0=tUv.
\end{equation}

Orthogonality follows from
\begin{equation}\label{eq:3.14}
UU^T=M^{-1}QQ^TM^{-1}=I,
\end{equation}
and hence $\|q_i\|=\alpha_i$, $\|q_0\|=t$. The aggregate is
\begin{equation}\label{eq:3.15}
\sum_iq_iu_i^T+q_0v^T=UQ^T=M^{-1}QQ^T=M.
\end{equation}

The scalar weights remain balanced, and
\begin{equation}\label{eq:3.16}
\det M=\sqrt{\det(QQ^T)}=|\det Q|=\det Q\ge R_n^{-1}.
\end{equation}

We now have every part of the certificate: normalization, balance, cone feasibility, positive definiteness, and a lower bound on its determinant.

\begin{proof}[Completion of the proof of Theorem~\ref{thm-sharp}]
Every ellipsoid in $K\cap H_v$ is feasible for the contact polytope and the cut. Lemma~\ref{lem-certificate} gives
\begin{equation}\label{eq:3.17}
\det A\le(\det M)^{-1}\le R_n.
\end{equation}

Since $w(K)=\kappa_n$ in John coordinates,
\begin{equation}\label{eq:3.18}
\frac{w(K\cap H_v)}{w(K)}
=\max_{\mathcal E(c,A)\subseteq K\cap H_v}\det A\le R_n.
\end{equation}

The ratio is affine invariant, proving \eqref{eq:3.1}, and the simplex model attains it. For $n=1$, write $K=[\ell,r]$; then $J(K)=K$, $c(K)=(\ell+r)/2$, and $w(K\cap H)=(r-\ell)/2=w(K)/2$. 
\end{proof}

The proof has also produced a quantitative object: the logarithmic determinant excess in \eqref{eq:3.12}. We next show that this excess is one part of an exact volume deficit. This is the starting point for equality and stability.

\section{Equality and stability}\label{sec:eq-stability}

With the sharp bound established, we examine the losses in its proof. The
exact deficit first gives the full equality family and then controls the
shape, center and weighted contacts of a nearly extremal ellipsoid. Throughout
this section we work in John coordinates in dimensions at least two.

\subsection{The exact deficit}\label{sec:4}

The contraction theorem gives an upper bound. To understand equality and near equality, we need to identify where the proof can lose volume. The certificate separates these losses into three terms.

\subsubsection{The two losses in the certificate inequality}\label{sec:4.1}

For a feasible primal--dual pair in Section~\ref{sec:2}, let
\begin{equation}\label{eq:4.1}
\tau_r=\beta_r-g_r^Tc,\qquad x_r=Ag_r,\qquad
\gamma_r=\alpha_r\tau_r-q_r^Tx_r\ge0,
\qquad Y=M^{1/2}AM^{1/2}.
\end{equation}

Each $\gamma_r$ is a cone pairing. The identities
$\operatorname{tr}Y=\operatorname{tr}(AM)$ and
$\det Y=\det A\det M$ give
\begin{equation}\label{eq:4.2}
\begin{aligned}
\mathcal C&:=\sum_r\gamma_r
=\sum_r\alpha_r\beta_r-c^T\sum_r\alpha_rg_r-\operatorname{tr}(AM)\\
&=n-\operatorname{tr}(AM),\\
\mathcal D&:=\operatorname{tr}Y-n-\log\det Y\ge0,\\
G&:=-\log\det M-\log\det A=\mathcal C+\mathcal D.
\end{aligned}
\end{equation}

Here $\mathcal C$ records the slack in the aggregated support inequalities. The matrix divergence $\mathcal D$ records the slack in the spectral inequality $\log\det Y\le\operatorname{tr}Y-n$. Its nonnegativity follows from $r-1-\log r\ge0$ for each positive eigenvalue $r$. The gap $G$ is computable from a feasible ellipsoid and a feasible certificate.

\subsubsection{Add the geometric loss}\label{sec:4.2}

For the explicit certificate of Section~\ref{sec:3}, define
$\Delta_{\rm geom}=\log(R_n\det M)=\log(R_n\det Q)$. Theorem~\ref{thm-frame-bound} proves that this term is nonnegative and controls $\chi$.

\begin{proposition}[{exact volume deficit}]\label{prop-gap}
In John coordinates, every feasible cut ellipsoid satisfies
\begin{equation}\label{eq:4.3}
\boxed{\epsilon:=\log\frac{R_n}{\det A}
=\Delta_{\rm geom}+\mathcal C+\mathcal D.}
\end{equation}

All three terms are nonnegative. For the actual new John ellipsoid, the same deficit is $\epsilon=\log(R_nw(K)/w(K\cap H))$ in the original coordinates.

\end{proposition}

\begin{proof}
Add the geometric loss to the primal--dual gap. The intermediate
certificate determinant cancels, leaving exactly the volume deficit:
\begin{equation}\label{eq:4.4}
\begin{aligned}
\Delta_{\rm geom}+\mathcal C+\mathcal D
&=\log(R_n\det M)-\log\det M-\log\det A\\
&=\log\frac{R_n}{\det A}=\epsilon.
\end{aligned}
\end{equation}

At the new optimum, $\det A=w(K\cap H_v)/w(K)$. Affine invariance gives the stated interpretation. 
\end{proof}

The three terms have distinct zero sets:
\begin{equation}\label{eq:4.5}
\begin{aligned}
\Delta_{\rm geom}=0&\iff\chi=0\iff s_i\in\{-1,b\}\quad\forall i,\\
\mathcal C=0&\iff\gamma_r=0\quad\forall r,\\
\mathcal D=0&\iff Y=I\iff A=M^{-1}.
\end{aligned}
\end{equation}

Thus zero deficit imposes the model contact geometry, complementarity, and shape matching simultaneously. For small deficit, the same decomposition gives quantitative control; in particular, $\chi\le1500\Delta_{\rm geom}\le1500\epsilon$.

\subsubsection{Small matrix divergence means close shape matching}\label{sec:4.3}

We use $\|\cdot\|_F$ for the Frobenius matrix norm; vector norms are Euclidean.

\begin{lemma}[{matrix matching}]\label{lem-matching}
For any normalized feasible pair, if $G\le1/4$, then
\begin{equation}\label{eq:4.6}
\|M^{1/2}AM^{1/2}-I\|_F\le2\sqrt G.
\end{equation}

\end{lemma}

\begin{proof}
We first use the small divergence to bound each eigenvalue from
above. On that bounded range, the same scalar divergence controls squared
distance from $1$. Put $h(t)=t-1-\log t$.
Integrating $h'(t)=1-1/t$ on either side of $1$ gives
\begin{equation}\label{eq:4.7}
h(t)\ge\frac{(t-1)^2}{2\max\{1,t\}}\quad(t>0),
\qquad h(2)=1-\log2>\frac14.
\end{equation}

For the eigenvalues $\rho_j>0$ of $Y$,
$\sum_jh(\rho_j)=\mathcal D\le G\le1/4$.
Since $h$ is increasing on $[1,\infty)$, each $\rho_j<2$. Hence
\begin{equation}\label{eq:4.8}
\|Y-I\|_F^2=\sum_j(\rho_j-1)^2
\le4\sum_jh(\rho_j)=4\mathcal D\le4G.
\end{equation}

Taking square roots proves the claim. 
\end{proof}

Section~\ref{sec:9} also gives an exact sum-of-squares identity for each cone pairing and a corresponding weighted estimate for multiplier directions. For the geometric conclusions below, the essential point is already visible: a small gap enforces matrix matching, while the geometric part of the deficit controls how close $M$ is to the model certificate.

\subsection{The equality family}\label{sec:5}

\subsubsection{The complete equality condition}\label{sec:5.1}

Continue in John coordinates and write
\begin{equation}\label{eq:5.1}
\mathcal Z_b=\{-v\}\cup\{u\in S^{n-1}:u^Tv=b\}.
\end{equation}

This set consists of the pole opposite the cut direction and a latitude sphere. The extremal ellipsoid $\mathcal E(c_*,A_*)$ is the one defined in \eqref{eq:3.2}.

\begin{theorem}[{equality classification}]\label{thm-equality}
A normalized body $K$ attains $w(K\cap H_v)=R_nw(K)$ if and only if it admits a John frame supported on $\mathcal Z_b$ such that
\begin{equation}\label{eq:5.2}
\operatorname{conv}\bigl(B_2^n\cup\mathcal E(c_*,A_*)\bigr)
\subseteq K\subseteq P_{\mathcal U}.
\end{equation}

The new John ellipsoid is then exactly $\mathcal E(c_*,A_*)$. If equality holds, every positive-weight John frame supporting $K$ lies in $\mathcal Z_b$. For a contact polytope $P_{\mathcal U}$, equality is equivalent simply to its defining John frame being supported on $\mathcal Z_b$.

\end{theorem}

\begin{proof}
For necessity, zero deficit will first determine the contact heights
and ellipsoid shape, then the active supports will determine its center.
For sufficiency, the candidate ellipsoid in the sandwich already supplies the
required volume.

\textbf{Necessity: recover the shape and then the center.} Fix any supporting John
frame and suppose $\det A=R_n$. The deficit identity gives
\begin{equation}\label{eq:5.3}
\begin{aligned}
\epsilon=0&\ \Longrightarrow\ \Delta_{\rm geom}=\mathcal C=\mathcal D=0\\
&\ \Longrightarrow\ \chi=0,\quad Q=M=M_0,\quad U=I,\quad A=A_*.
\end{aligned}
\end{equation}

The shape is now fixed. To recover the center, compare the active support
values with those of the model. Since $\alpha_i>0$ and
$p_i=A_*u_i/\|A_*u_i\|$, complementarity reduces to
\begin{equation}\label{eq:5.4}
0=\gamma_i
=\alpha_i\bigl(1-u_i^Tc-\|A_*u_i\|\bigr),
\qquad
u_i^T(c-c_*)=0\quad\forall i.
\end{equation}

The second identity uses the model support equations. Isotropy now recovers the center:
\begin{equation}\label{eq:5.5}
\|c-c_*\|^2
=\sum_i\lambda_i\bigl(u_i^T(c-c_*)\bigr)^2=0.
\end{equation}

Thus $\mathcal E(c_*,A_*)\subseteq K\cap H_v$, and convexity yields the left inclusion in \eqref{eq:5.2}. The right inclusion is \eqref{eq:2.18}. The frame was arbitrary, proving the asserted quantifier over every John frame.

\textbf{Sufficiency: use the included model ellipsoid.} Suppose \eqref{eq:5.2} holds for a
John frame supported on $\mathcal Z_b$. John's characterization and the model
feasibility give
\begin{equation}\label{eq:5.6}
w(K)=\kappa_n,\qquad
w(K\cap H_v)\ge\kappa_n\det A_*=R_n\kappa_n.
\end{equation}

Theorem~\ref{thm-sharp} supplies the reverse inequality. For $K=P_{\mathcal U}$, model feasibility follows directly from the three support values in Section~\ref{sec:3.2}. 
\end{proof}

The equality condition determines the new ellipsoid and the active contact geometry. It leaves freedom in the rest of the body, as the two inclusions in \eqref{eq:5.2} show.

\subsubsection{Extremizers beyond regular simplices}\label{sec:5.2}

The latitude part of an equality frame can be built from any John frame in one lower dimension. This gives a concrete description of the freedom in Theorem~\ref{thm-equality}.

\begin{proposition}[{transverse construction}]\label{prop-extremizers}
Let unit vectors $w_i\in v^\perp$ and positive weights $\mu_i$ satisfy
$\sum_i\mu_iw_i=0$ and $\sum_i\mu_iw_iw_i^T=P$. Then the contacts
\begin{equation}\label{eq:5.7}
u_i=bv+\sqrt{1-b^2}\,w_i,\qquad
\lambda_i=\frac{\mu_i}{1-b^2},\qquad
u_{\rm pole}=-v,\qquad \lambda_{\rm pole}=\frac1{1+b}
\end{equation}
form a John frame supported on $\mathcal Z_b$.

\end{proposition}

\begin{proof}
The weight rescaling compensates for the transverse length
$\sqrt{1-b^2}$, while the pole weight balances the axial first moment.
We check that the same choices also give the axial second moment.
The vectors have unit norm because $w_i\perp v$ and
$\|u_i\|^2=b^2+(1-b^2)\|w_i\|^2=1$. Taking traces in the transverse second moment gives
\begin{equation}\label{eq:5.8}
\sum_i\mu_i=n-1=\frac{1-b}{b},\qquad
S:=\sum_i\lambda_i=\frac1{b(1+b)}.
\end{equation}

The pole and latitude contributions therefore satisfy
\begin{equation}\label{eq:5.9}
\begin{aligned}
\lambda_{\rm pole}u_{\rm pole}+\sum_i\lambda_i u_i
&=\left(-\frac1{1+b}+bS\right)v
+\frac1{\sqrt{1-b^2}}\sum_i\mu_iw_i=0,\\
\lambda_{\rm pole}u_{\rm pole}u_{\rm pole}^T+\sum_i\lambda_i u_i u_i^T
&=\left(\frac1{1+b}+b^2S\right)vv^T+\sum_i\mu_iw_iw_i^T\\
&=vv^T+P=I.
\end{aligned}
\end{equation}

The mixed terms in the second line vanish by $\sum_i\mu_iw_i=0$. All weights are positive and all heights lie in $\{-1,b\}$. 
\end{proof}

For example, take $w_i=\pm e_j$, $j=1,\ldots,n-1$, with $\mu_i=1/2$. In dimension three, set $v=e_3$. The four latitude contacts have John weight $9/16$ each, and the pole $-e_3$ has weight $3/4$. Their contact polytope is the square pyramid
\begin{equation}\label{eq:5.10}
K=\{(x,y,z):z\ge-1,\quad
2\sqrt2|x|+z\le3,\quad 2\sqrt2|y|+z\le3\}.
\end{equation}

Its base vertices are $(\pm\sqrt2,\pm\sqrt2,-1)$ and its apex is $(0,0,3)$. Proposition~\ref{prop-extremizers} and Theorem~\ref{thm-equality} give $J(K)=B_2^3$ and
\begin{equation}\label{eq:5.11}
\begin{gathered}
J(K\cap\{z\le0\})=c_*+A_*B_2^3,\qquad
c_*=-\frac12e_3,\\
A_*=\operatorname{diag}\!\left(\sqrt{\frac32},\sqrt{\frac32},\frac12\right),
\qquad \frac{w(K\cap\{z\le0\})}{w(K)}=\det A_*=\frac34=R_3.
\end{gathered}
\end{equation}

\begin{figure}[tbp]
\centering
\includegraphics[width=0.98\linewidth]{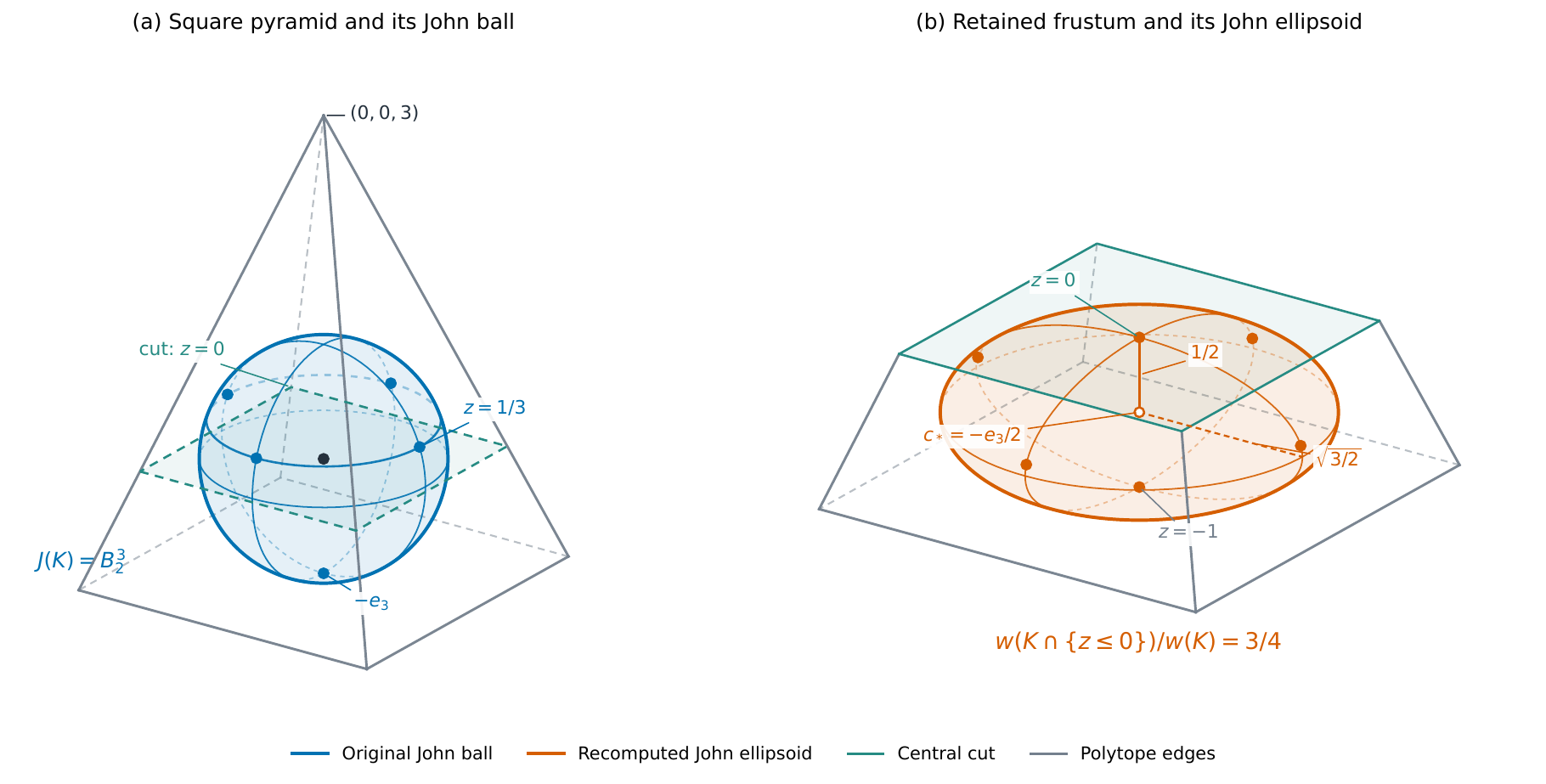}
\caption{A square-pyramid extremizer in dimension three. (a) The unit ball touches four triangular side facets at height $z=1/3$ and the square base at $-e_3$. (b) The central cut leaves a square frustum. Its John ellipsoid touches the four sides, the base, and the cut; its semiaxes are $\sqrt{3/2},\sqrt{3/2},1/2$. Filled colored points mark contacts, and the open orange point marks the new center. Both panels use the same orthographic viewing direction, with an enlarged view in (b). This five-facet extremizer exhibits the freedom beyond the simplex model in Figure~\ref{fig:1}.}\label{fig:2}
\end{figure}

The regular tetrahedron is therefore only one member of the three-dimensional equality family. Any body between the two sets in \eqref{eq:5.2} is also an extremizer.

This has a direct consequence for stability. Near equality should force proximity to the common extremal ellipsoid and contact family. It cannot force every body to be close to a unique simplex.

\subsection{Stability of the ellipsoid and its contacts}\label{sec:6}

The deficit identity provides the qualitative equality argument almost immediately. Quantitative stability requires estimates that turn its three small terms into distances. We state the geometric estimates first, then combine them with matrix matching.

\subsubsection{The geometric estimates needed by the proof}\label{sec:6.1}

\begin{lemma}[{geometric control}]\label{lem-geometric-control}
For the explicit certificate,
\begin{equation}\label{eq:6.1}
\|Q-M_0\|_F^2\le\frac{1-b}{2-b}\chi,\qquad
\sum_i\lambda_i\operatorname{dist}(u_i,\mathcal Z_b)^2
\le\frac{32}{3}\chi.
\end{equation}

If $\chi\le3/2000$, then $\Delta_{\rm geom}\ge\chi/32$.

\end{lemma}

\textbf{Proof deferred to Sections~\ref{sec:9.2}--\ref{sec:9.3}.} Section~\ref{sec:9.2} proves the two global
estimates through Lemma~\ref{lem-residual}. Section~\ref{sec:9.3} proves the local implication through
Lemma~\ref{lem-local}, checks the conversion between the two smallness conditions, and
completes Lemma~\ref{lem-geometric-control}.

Even without the small-defect assumption, \eqref{eq:4.3} and Lemma~\ref{lem-geometric-control} give
$\sum_i\lambda_i\operatorname{dist}(u_i,\mathcal Z_b)^2\le16000\epsilon$.
Thus the contact geometry is controlled globally, with a better constant near equality.

\subsubsection{Shape, center, and contact stability}\label{sec:6.2}

\begin{theorem}[{dimension-independent stability}]\label{thm-stability}
Let $\mathcal E(c,A)$ be any feasible cut ellipsoid in John coordinates, and assume
$0\le\epsilon=\log(R_n/\det A)\le10^{-6}$. Then
\begin{equation}\label{eq:6.2}
\begin{aligned}
\chi&\le32\epsilon,&
\|Q-M_0\|_F&\le4\sqrt\epsilon,\\
\|A-A_*\|_F&\le30\sqrt\epsilon,&
\|c-c_*\|&\le50\sqrt\epsilon,
\end{aligned}
\qquad
\sum_i\lambda_i\operatorname{dist}(u_i,\mathcal Z_b)^2
\le\frac{1024}{3}\epsilon.
\end{equation}

\end{theorem}

The constants are independent of $n$. The threshold $10^{-6}$ is a sufficient
condition for entering the local regime of Lemma~\ref{lem-geometric-control}; the proof below shows
this implication explicitly. The numbers $30$ and $50$ are convenient upper
bounds obtained when the shape and center estimates are combined.

The proof below uses Lemma~\ref{lem-geometric-control}, whose proof is deferred as stated above.
The polar-factor and support-mismatch estimates used in the shape and center
steps are established in Section~\ref{sec:9.4}.

\begin{proof}
We pass from the contact defect to the certificate, then to the shape and center.

\textbf{Contact geometry.} First use the global bound to enter the small-defect
regime, then use the stronger local bound. This is the purpose of the
threshold on $\epsilon$. By \eqref{eq:4.3} and Lemma~\ref{lem-geometric-control},
\begin{equation}\label{eq:6.3}
\chi\le1500\epsilon\le\frac3{2000}
\quad\Longrightarrow\quad
\chi\le32\Delta_{\rm geom}\le32\epsilon.
\end{equation}

Consequently,
\begin{equation}\label{eq:6.4}
\|Q-M_0\|_F^2\le\frac{1-b}{2-b}\chi\le\frac\chi2\le16\epsilon,
\qquad
\sum_i\lambda_i\operatorname{dist}(u_i,\mathcal Z_b)^2
\le\frac{32}{3}\chi\le\frac{1024}{3}\epsilon.
\end{equation}

\textbf{Shape.} There are two errors to combine: matrix matching controls $A$ relative
to $M^{-1}$, and the geometric estimate controls $M^{-1}$ relative to
$A_*=M_0^{-1}$. Set $Y=M^{1/2}AM^{1/2}$.
Then $G\le\epsilon$ and Lemma~\ref{lem-matching} gives $\|Y-I\|_F\le2\sqrt\epsilon$. Section~\ref{sec:9.4} establishes
$\lambda_{\min}(M),\lambda_{\min}(M_0)>7/10$ and
$\|M-M_0\|_F\le12\sqrt\epsilon$. Using $A_*=M_0^{-1}$,
\begin{equation}\label{eq:6.5}
\begin{aligned}
\|A-A_*\|_F
&\le\|M^{-1/2}(Y-I)M^{-1/2}\|_F
+\|M^{-1}(M_0-M)M_0^{-1}\|_F\\
&\le\frac{\|Y-I\|_F}{\lambda_{\min}(M)}
+\frac{\|M-M_0\|_F}{\lambda_{\min}(M)\lambda_{\min}(M_0)}\\
&\le\left(\frac{20}{7}+\frac{1200}{49}\right)\sqrt\epsilon
\le30\sqrt\epsilon.
\end{aligned}
\end{equation}

The inverse-difference identity preserves the multiplication order and requires no commutativity.

\textbf{Center.} Recover the displacement from its projections on the contact
normals. Each projection has three contributions: support slack, mismatch
between the contact height and the equality model, and the shape error already
bounded above. Let $r_i=1-c^Tu_i-\|Au_i\|\ge0$. Since
$\gamma_i\ge\alpha_i r_i$ and $\alpha_i\ge\lambda_i/4$,
\begin{equation}\label{eq:6.6}
\sum_i\lambda_i r_i\le4\sum_i\alpha_i r_i\le4\mathcal C\le4\epsilon.
\end{equation}

The function $g_*(s)=\|A_*u\|-1-s/2$, for $u^Tv=s$, measures the model's
support mismatch at that contact height. It vanishes at the equality contacts.
Section~\ref{sec:9.4} proves
$\sum_i\lambda_i g_*(s_i)^2\le8\chi$. Subtracting support values and using isotropy gives
\begin{equation}\label{eq:6.7}
\begin{aligned}
c-c_*={}&-\sum_i\lambda_i r_i u_i
-\sum_i\lambda_i g_*(s_i)u_i\\
&+\sum_i\lambda_i\bigl(\|A_*u_i\|-\|Au_i\|\bigr)u_i.
\end{aligned}
\end{equation}

For any scalars $a_i$, the frame Bessel inequality is
$\|\sum_i\lambda_i a_i u_i\|^2\le\sum_i\lambda_i a_i^2$.
Moreover,
\begin{equation}\label{eq:6.8}
\sum_i\lambda_i\bigl|\|A_*u_i\|-\|Au_i\|\bigr|^2
\le\sum_i\lambda_i\|(A_*-A)u_i\|^2=\|A_*-A\|_F^2.
\end{equation}

The triangle and Bessel inequalities now yield
\begin{equation}\label{eq:6.9}
\|c-c_*\|\le4\epsilon+\|A-A_*\|_F+\sqrt{8\chi}
\le50\sqrt\epsilon.
\end{equation}

This proves all estimates in \eqref{eq:6.2}. 
\end{proof}

The conclusion applies to every nearly extremal feasible ellipsoid, including the optimizer. Its geometric interpretation is that a small logarithmic volume deficit forces a square-root scale error in the shape and center. The contact estimate is weighted, so it does not control a single contact whose John weight is arbitrarily small. All distances are measured after John normalization, and no whole-body distance to a unique extremizer is asserted.

\subsection{Supporting estimates for stability}\label{sec:9}

We now supply the proof of Lemma~\ref{lem-geometric-control} deferred in Section~\ref{sec:6.1}, and the
auxiliary shape and center estimates deferred in the proof of Theorem~\ref{thm-stability}.
Sections~\ref{sec:9.2}--\ref{sec:9.3} prove the global and local parts of Lemma~\ref{lem-geometric-control}, respectively;
Section~\ref{sec:9.4} supplies the remaining inputs to Theorem~\ref{thm-stability}. The quantities have
three distinct roles: $\epsilon$ measures the total logarithmic volume deficit,
$\chi$ measures departure of the contacts from the equality set, and
$\sigma=2\chi/[d(2+b)]$ is the same contact defect in the block normalization.
We first turn $\chi$ into matrix and contact-distance bounds. We then improve
the determinant lower bound near equality and transfer the resulting matrix
control to the ellipsoid's shape and center.

The constants record successive sufficient estimates in this chain.
The threshold in Lemma~\ref{lem-local} ensures that a quadratic logarithm bound applies;
Section~\ref{sec:9.4} uses fixed spectral bounds to control inversion. The final
coefficients $30$ and $50$ collect those estimates and are not asserted to
be optimal stability constants.

\subsubsection{Complementarity and matrix matching}\label{sec:9.1}

Use the primal and dual quantities from \eqref{eq:4.1}--\eqref{eq:4.2}. The following identities
also give the weighted multiplier estimate mentioned in Section~\ref{sec:4.3}.
Here $(\tau,x)$ denotes a generic scalar--vector cone pair; the indexed
$(\tau_r,x_r)$ are the support slack and vector from \eqref{eq:4.1}.

\begin{proposition}[{quantitative complementarity and matrix matching}]\label{prop-gap-control}
For cone-feasible $(\alpha,q)$ and $(\tau,x)$,
\begin{equation}\label{eq:9.1}
2\alpha\tau(\alpha\tau-q^Tx)
=\|\alpha x-\tau q\|^2
+\alpha^2(\tau^2-\|x\|^2)+\tau^2(\alpha^2-\|q\|^2).
\end{equation}

If $G\le\epsilon\le1/4$, then
\begin{equation}\label{eq:9.2}
\|Y-I\|_F\le2\sqrt\epsilon,\qquad
\sum_{\alpha_r>0}\alpha_r\tau_r
\left\|\frac{x_r}{\tau_r}-\frac{q_r}{\alpha_r}\right\|^2
\le2\epsilon.
\end{equation}

Here zero-normal redundant inequalities are removed. The first estimate also holds with $Y$ replaced by $A^{1/2}MA^{1/2}$.

\end{proposition}

\begin{proof}
The sum-of-squares identity converts a small cone pairing into
closeness of the normalized cone directions. Expanding the squares proves
\eqref{eq:9.1}, and cone feasibility makes every term on its right nonnegative.
For a nonzero normal, $A\succ0$ and primal feasibility imply
$\tau_r\ge\|Ag_r\|>0$. Divide \eqref{eq:9.1} by $\alpha_r\tau_r$ when $\alpha_r>0$
and sum. The total cone gap is at most $\epsilon$, giving the second bound.

The matrix estimate is Lemma~\ref{lem-matching}:
$\|Y-I\|_F\le2\sqrt G\le2\sqrt\epsilon$.
The symmetric positive-definite matrices $M^{1/2}AM^{1/2}$ and
$A^{1/2}MA^{1/2}$ have the same eigenvalues, so their Frobenius distances
from the identity agree. This proves the final assertion. 
\end{proof}

\subsubsection{Global part of Lemma~\ref{lem-geometric-control}: geometric residuals}\label{sec:9.2}

We return to the block and contact notation of Section~\ref{sec:8}. The aim is to make
the scalar defect $\chi$ useful in two ways: as a bound on the certificate's
matrix error, and as a bound on the distance of the contacts from the equality
set. For the matrix bound, we seek a residual for each contact whose squared
norm is proportional to that contact's summand in $\chi$. Isotropy will then
combine these residuals without a dimension factor.

\begin{lemma}[{geometric residual}]\label{lem-residual}
For the explicit certificate,
\begin{equation}\label{eq:9.3}
\|Q-M_0\|_F^2\le\frac{1-b}{d}\chi,\qquad
\|\widehat F-I\|_F^2\le\frac\chi d,
\qquad
\sum_i\lambda_i\operatorname{dist}(u_i,\mathcal Z_b)^2
\le\frac{32}{3}\chi.
\end{equation}

\end{lemma}

\begin{proof}
\textbf{Step 1: express the certificate error as contact residuals.} Subtract the model certificate from $Q$ and use the zero first moment to add
a term that vanishes in the aggregate. This makes the individual residual
vanish at both the pole and the model latitude.
For $u=sv+y$ with $y\perp v$, the resulting residual is
\begin{equation}\label{eq:9.4}
\xi(u)=-\frac{s-b}{d(1-as)}
\bigl[(1-b)(s+1)v+\sqrt{1-b}\,y\bigr].
\end{equation}

Expansion of $Q-M_0$, followed by the first and second John moments, gives
\begin{equation}\label{eq:9.5}
Q-M_0=\sum_i\lambda_i u_i\xi(u_i)^T,
\qquad
\|\xi(u)\|^2=\frac{1-b}{d}\frac{(1+s)(s-b)^2}{1-as}.
\end{equation}

For clarity, the axial coefficient before simplification is
$f(s)(s-a)-(2-t)s+b$, and the transverse coefficient is
$kf(s)-\sqrt{1-b}$. The added constant vector $bv$ contributes zero to the aggregate by balance. Their simplification gives \eqref{eq:9.4}, and orthogonality gives the norm identity.

\textbf{Step 2: use isotropy to bound the total matrix error.} The matrix with
columns $\sqrt{\lambda_i}u_i$ has orthonormal rows. Its Bessel inequality
bounds the squared Frobenius norm of the aggregate by the sum of the
weighted squared residuals. Applying this to \eqref{eq:9.5} proves the first estimate.
To transfer it to the normalized matrix, use
$Q-M_0=F_0^{1/2}(\widehat F-I)F_0^{1/2}D$.
The product of the least singular values of the two outside factors is
$\sqrt{1-b}$, giving the second estimate. The two matrix constants in \eqref{eq:9.3}
thus come from the exact residual identity and the fixed normalization.

\textbf{Step 3: compare the scalar residual with geometric distance.}
The equality set consists of a pole and a latitude sphere. We choose an
explicit comparison point in one of these two pieces. It suffices to show
\begin{equation}\label{eq:9.6}
\frac{(1+s)(s-b)^2}{1-as}
\ge\frac3{32}\operatorname{dist}(u,\mathcal Z_b)^2.
\end{equation}

The split at $s=-1/2$ gives simple uniform bounds in both regions.
If $s\le-1/2$, use the pole, whose squared distance is $2(1+s)$;
then $(s-b)^2\ge1/4$ and $1-as\le4/3$.
If $s\ge-1/2$, choose the point on the latitude with the same transverse
direction, choosing any unit transverse vector if $y=0$. Its squared distance is
\begin{equation}\label{eq:9.7}
(s-b)^2+\bigl(\sqrt{1-s^2}-\sqrt{1-b^2}\bigr)^2\le4(s-b)^2.
\end{equation}

Indeed, rationalization bounds the squared root difference by
$(s-b)^2(s+b)^2/(1-b^2)\le3(s-b)^2$.
Also $(1+s)/(1-as)\ge3/8$. These two comparisons give the same coefficient
$3/32$ in \eqref{eq:9.6}; its reciprocal is the contact-distance constant in \eqref{eq:9.3}.
Summing with the John weights finishes the proof. 
\end{proof}

\textbf{Global part of Lemma~\ref{lem-geometric-control} proved.} The first and third inequalities in
\eqref{eq:9.3} are exactly the two estimates in \eqref{eq:6.1}, since $d=2-b$.
The local implication in Lemma~\ref{lem-geometric-control} remains to be proved in Section~\ref{sec:9.3}.

\textbf{What the global estimates provide.} Combining Theorem~\ref{thm-determinant}, Proposition~\ref{prop-gap},
and Lemma~\ref{lem-residual} gives
\begin{equation}\label{eq:9.8}
\chi\le1500\epsilon,\qquad
\sum_i\lambda_i\operatorname{dist}(u_i,\mathcal Z_b)^2\le16000\epsilon,
\qquad
\log\det\widehat F\ge\frac1{1000}\|\widehat F-I\|_F^2.
\end{equation}

\subsubsection{Local part of Lemma~\ref{lem-geometric-control}: coercivity near equality}\label{sec:9.3}

The global estimate is enough to place a nearly extremal ellipsoid in a small
neighborhood of the model. Here we improve its coefficient within that
neighborhood. The idea is that the positive linear term
$\operatorname{tr}(\widehat F-I)=\sigma-x$ should exceed the quadratic remainder
in $\log\det\widehat F$. We need to bound the whole perturbation, including
$j$, tightly enough for this difference to remain positive.

The threshold $10^{-3}$ below is a sufficient size condition for that
comparison. The intermediate bounds $1/40$ and $1/30$ specify the small
spectral intervals used in the proof; the useful conclusion is the final
lower bound $\sigma/16$.

\begin{lemma}[{local coercivity}]\label{lem-local}
If $\sigma\le10^{-3}$, then
\begin{equation}\label{eq:9.9}
\log\det\widehat F\ge\frac\sigma{16}\ge\frac\chi{32}.
\end{equation}

\end{lemma}

\begin{proof}
If $\sigma=0$, the assertion follows from $\widehat F=I$.
Suppose $\sigma>0$.

\textbf{Step 1: confine the transverse eigenvalues.} We first use Lemma~\ref{lem-variational} to
control $E$ and $j$ in the metric determined by $K_E$.
Here $0\le L\le d(1-b)/(d+b^2)<1$ by its definition and
$(d+b^2)-d(1-b)=2b>0$. Lemma~\ref{lem-variational} therefore gives
\begin{equation}\label{eq:9.10}
A_E+\|K_E^{-1}j\|^2\le\frac\sigma2.
\end{equation}

Every eigenvalue $e>-1/2$ of $E$ therefore satisfies
$e^2/(1+2e)\le1/2000$. This scalar inequality implies $|e|\le1/40$:
the quotient increases as $e$ moves away from zero on either side, and its
values at $\pm1/40$ already exceed $1/2000$.
The choice $1/40$ gives the simple interval $[19/20,21/20]$ for the eigenvalues
of $K_E$. We can now compare the metric in \eqref{eq:9.10} with ordinary Frobenius
and Euclidean norms:
\begin{equation}\label{eq:9.11}
\|E\|_F^2+\frac{20}{21}\|j\|^2\le\frac{21}{40}\sigma.
\end{equation}

\textbf{Step 2: control the full perturbation.} Its squared norm contains
$2\|j\|^2$, so the estimate in \eqref{eq:9.11} needs one further bound on the coupling.
The Gram constraint gives $\|j\|^2\le x(t+1/40)\le47\sigma/160$.
The coefficient $47/160$ combines $xt\le9\sigma/32$ and
$x/40\le\sigma/80$. Also $x^2\le\sigma/4000$.
Combining these estimates leaves the convenient upper bound $5\sigma/6$:
\begin{equation}\label{eq:9.12}
\|\widehat F-I\|_F^2=x^2+2\|j\|^2+\|E\|_F^2\le\frac56\sigma,
\qquad \|\widehat F-I\|_{\rm op}<\frac1{30}.
\end{equation}

\textbf{Step 3: let the trace gain dominate the logarithmic remainder.}
The size assumption and \eqref{eq:9.12} put every eigenvalue perturbation in
$[-1/30,1/30]$. On this interval, integration of the derivative gives
$\log(1+z)\ge z-(15/29)z^2$; the coefficient is
$15/29=1/[2(1-1/30)]$.
The quadratic term is controlled by \eqref{eq:9.12}, while
$\operatorname{tr}(\widehat F-I)=(1-d/4)\sigma\ge\sigma/2$.
Their difference is
\begin{equation}\label{eq:9.13}
\log\det\widehat F\ge\frac\sigma2-\frac{15}{29}\frac{5\sigma}{6}
=\frac{2\sigma}{29}\ge\frac\sigma{16}.
\end{equation}

The bound $2\sigma/29$ is rounded down to $\sigma/16$ for subsequent use.
Finally, $\sigma=2\chi/(4-b^2)\ge\chi/2$. 
\end{proof}

\textbf{Completion of the deferred proof of Lemma~\ref{lem-geometric-control}.} The first two bounds in Lemma~\ref{lem-geometric-control}
are part of Lemma~\ref{lem-residual}. Its hypothesis is expressed in $\chi$, while the local
lemma uses $\sigma$, so we check the conversion explicitly. For the small-defect
assertion, $\chi\le3/2000$ implies
$\sigma=2\chi/(4-b^2)\le(8/15)\chi\le1/1250<10^{-3}$.
Lemma~\ref{lem-local} and \eqref{eq:8.2} therefore give $\Delta_{\rm geom}\ge\chi/32$.
This proves the remaining assertion and completes Lemma~\ref{lem-geometric-control}, as promised in Section~\ref{sec:6.1}.

\subsubsection{Remaining estimates for Theorem~\ref{thm-stability}: shape and center}\label{sec:9.4}

We now prove the polar-factor and support-mismatch estimates deferred in
Section~\ref{sec:6.2}. We finish the passage from certificate error to ellipsoid error under the
hypotheses of Theorem~\ref{thm-stability}, in particular $\epsilon\le10^{-6}$.
The contact estimate \eqref{eq:6.3} already gives $\chi\le32\epsilon$ and hence
$\|Q-M_0\|_F\le4\sqrt\epsilon$. What remains is to transfer this bound through
the polar factor and matrix inversion, then recover the center from support
values. The constants below are explicit bounds for these transfers.

\textbf{Step 1: keep the polar factor uniformly invertible.} Set
$\delta=\|Q-M_0\|_{\rm op}\le4\sqrt\epsilon\le0.004$. Since
$M_0=2vv^T+\sqrt{1-b}P$, the triangle inequality gives, for every vector $x$,
\begin{equation}\label{eq:9.14}
\left(\frac1{\sqrt2}-\delta\right)\|x\|
\le\|Q^Tx\|\le(2+\delta)\|x\|,
\qquad
\|Q^Tx\|^2=x^TM^2x.
\end{equation}

Apply these inequalities to unit eigenvectors of $M\succ0$. Every eigenvalue
lies in $[1/\sqrt2-\delta,\,2+\delta]$.
The rational bounds $7/10$ and $201/100$ lie outside this interval uniformly
under the assumed bound on $\delta$:
\begin{equation}\label{eq:9.15}
\lambda_{\min}(M)>\frac7{10},\qquad
\lambda_{\max}(M)<\frac{201}{100},\qquad
\|Q\|_{\rm op}\le2+\delta<\frac{201}{100}.
\end{equation}

The lower bound prevents matrix inversion from magnifying small errors
without control. The upper bound will control the product difference in the
next step.

\textbf{Step 2: transfer the error from $Q$ to its polar factor.} Since $M^2=QQ^T$,
the difference $X=M-M_0$ satisfies a linear matrix equation. Taking its
Frobenius inner product with $X$ uses the positive spectral lower bounds to
control $\|X\|_F$. The identity $MX+XM_0=QQ^T-M_0^2$ gives
\begin{equation}\label{eq:9.16}
\begin{aligned}
(\lambda_{\min}(M)+\lambda_{\min}(M_0))\|X\|_F^2
&\le\operatorname{tr}(X^TMX)+\operatorname{tr}(X^TXM_0)\\
&=\langle X,QQ^T-M_0^2\rangle_F\\
&\le\|X\|_F\|QQ^T-M_0^2\|_F.
\end{aligned}
\end{equation}

For $X=0$ the conclusion is immediate; otherwise cancel $\|X\|_F$. Since
$QQ^T-M_0^2=(Q-M_0)Q^T+M_0(Q^T-M_0)$,
\begin{equation}\label{eq:9.17}
\begin{aligned}
\|M-M_0\|_F
&\le\frac{\|QQ^T-M_0^2\|_F}
{\lambda_{\min}(M)+\lambda_{\min}(M_0)}\\
&\le\frac{(\|Q\|_{\rm op}+\|M_0\|_{\rm op})\|Q-M_0\|_F}{7/5}
\le\frac{401}{35}\sqrt\epsilon
\le12\sqrt\epsilon.
\end{aligned}
\end{equation}

Here $401/35$ is the coefficient obtained by inserting the spectral bounds;
$12$ is a convenient upper bound for it.

\textbf{Step 3: transfer the error to the ellipsoid shape.} Matrix matching
controls $A-M^{-1}$, and Step 2 controls $M^{-1}-M_0^{-1}$ using the same
positive lower bounds. The inverse-difference estimate in \eqref{eq:6.5} gives
\begin{equation}\label{eq:9.18}
\|A-A_*\|_F
\le\left(\frac{20}{7}+\frac{1200}{49}\right)\sqrt\epsilon
\le30\sqrt\epsilon.
\end{equation}

The two terms in the coefficient correspond to these two errors, and their
sum is less than $30$. The non-strict conclusion includes $\epsilon=0$.

\textbf{Step 4: express the model's support mismatch through the contact defect.}
To recover the center, we need to compare the model support value with the
original supporting plane at each contact. Squaring and subtracting exposes
the pole and latitude factors directly. For a unit vector $u=sv+y$, $y\perp v$,
compute
\begin{equation}\label{eq:9.19}
\|A_*u\|^2-\left(1+\frac s2\right)^2
=\frac{s^2}{4}+\frac{1-s^2}{1-b}-\left(1+\frac s2\right)^2
=\frac{(1+s)(b-s)}{1-b}.
\end{equation}

Rationalizing therefore gives
\begin{equation}\label{eq:9.20}
g_*(s):=\|A_*u\|-1-\frac s2
=\frac{(1+s)(b-s)}{(1-b)(\|A_*u\|+1+s/2)}.
\end{equation}

The squared mismatch is therefore the contact-defect summand times a bounded
factor. The bounds $\|A_*u\|+1+s/2\ge1$, $(1+s)(1-as)\le2$, and
$(1-b)^2\ge1/4$ give the coefficient $8$:
\begin{equation}\label{eq:9.21}
g_*(s)^2
\le8\,\frac{(1+s)(s-b)^2}{1-as},
\qquad
\sum_i\lambda_i g_*(s_i)^2\le8\chi.
\end{equation}

\textbf{Step 5: recover the center and collect the estimates.} Subtracting support
values expresses each projection of $c-c_*$ as the sum of the slack, the
model mismatch, and the shape error:
\begin{equation}\label{eq:9.22}
(c-c_*)^Tu_i=-r_i-g_*(s_i)+\|A_*u_i\|-\|Au_i\|.
\end{equation}

Multiplication by $\lambda_i u_i$ and summation uses isotropy to recover
$c-c_*$ itself. The triangle and Bessel estimates in Section~\ref{sec:6} give
\begin{equation}\label{eq:9.23}
\|c-c_*\|\le4\epsilon+\sqrt{8\chi}+\|A-A_*\|_F
\le50\sqrt\epsilon.
\end{equation}

In this last estimate, $\sqrt{8\chi}\le16\sqrt\epsilon$, the shape term is
at most $30\sqrt\epsilon$, and $4\epsilon\le0.004\sqrt\epsilon$.
The coefficient $50$ is a single convenient upper bound for their sum.
Together with the completed proof of Lemma~\ref{lem-geometric-control}, these estimates justify
all deferred inputs in the proof of Theorem~\ref{thm-stability}.

\section{Exact duality and further consequences}\label{sec:further}

The geometric proof used one feasible certificate. We now establish the
exact primal--dual interpretation and examine its computational and matrix
consequences. Exact duality and iteration apply in every positive dimension;
the subsequent results for the explicit frame certificate use the definitions
from Section~\ref{sec:3} in dimensions at least two.

\subsection{Exact certificates and optimal contacts}\label{sec:2.3}

Optimizing the finite certificate recovers the ellipsoid optimum. The
construction below uses John contacts at an actual maximizing ellipsoid and
also identifies the contacts recovered from any optimal primal--dual pair.

\textbf{Construct an optimal certificate.} Return to the bounded full-dimensional
polytope $\mathcal P$ of \eqref{eq:2.3}. The next proposition proves attainment and
equality of the primal and normalized dual values.

\begin{proposition}[{optimal multipliers from John contacts}]\label{prop-exact-duality}
Problems \eqref{eq:2.3} and \eqref{eq:2.8} attain their optima. If $\mathcal E(c,A)=J(\mathcal P)$, there are normalized dual multipliers with aggregate $M=A^{-1}$ and zero cone pairings. Hence the optimal values of \eqref{eq:2.3} and \eqref{eq:2.8} sum to zero.

\end{proposition}

\begin{proof}
Construct the dual variables from the contacts of the actual
optimal ellipsoid. We first identify the active facet normal for each contact,
then rescale the John weights so that their moments become the dual constraints.

\textbf{Identify the active facets.} Let $\mathcal E(c,A)=J(\mathcal P)$,
whose existence is classical, and apply John's characterization to
$A^{-1}(\mathcal P-c)$. A unit contact $p$ lies on some defining plane.
For its index $r$,
\begin{equation}\label{eq:2.10}
\|Ag_r\|\le\beta_r-g_r^Tc=(Ag_r)^Tp
\le\|Ag_r\|\|p\|=\|Ag_r\|.
\end{equation}

Hence $p=Ag_r/\|Ag_r\|$ and the plane is active. Assign contacts to such indices and combine their weights. Writing $\ell_r=\|Ag_r\|>0$, this gives
\begin{equation}\label{eq:2.11}
p_r=\frac{Ag_r}{\ell_r},\qquad
\sum_r\mu_rp_r=0,\qquad
\sum_r\mu_rp_rp_r^T=I,\qquad \sum_r\mu_r=n.
\end{equation}

\textbf{Turn contact weights into multipliers.} The weights satisfy $\mu_r\ge0$ and
$\mu_r(\beta_r-g_r^Tc-\ell_r)=0$. Dividing by the active support length
$\ell_r$ makes the zero first moment become dual balance.
Define $\alpha_r=\mu_r/\ell_r$ and $q_r=\alpha_rp_r$. Then
\begin{equation}\label{eq:2.12}
\begin{aligned}
\sum_r\alpha_rg_r&=A^{-1}\sum_r\mu_rp_r=0,\\
\sum_r\alpha_r\beta_r
&=c^T\sum_r\alpha_rg_r+\sum_r\mu_r=n,\\
\|q_r\|&=\alpha_r,\\
\alpha_r(\beta_r-g_r^Tc)-q_r^TAg_r
&=\alpha_r(\beta_r-g_r^Tc-\ell_r)=0.
\end{aligned}
\end{equation}

\textbf{Use isotropy to match the inverse shape.} The second moment supplies the
remaining dual matrix equation. Since $g_r=\ell_rA^{-1}p_r$, the aggregate is
\begin{equation}\label{eq:2.13}
\sum_rq_rg_r^T=\sum_r\mu_rp_rp_r^TA^{-1}=A^{-1}.
\end{equation}

The cone constraints and zero pairings follow immediately. Lemma~\ref{lem-certificate} now proves that this normalized certificate is dual optimal and that the gap is zero. 
\end{proof}

\textbf{Recover contacts from an optimal pair.} The converse identifies the active
supports and shows that the normalized directions recover a John frame.

\begin{proposition}[{optimality and contact recovery}]\label{prop-kkt}
A feasible pair for \eqref{eq:2.3} and \eqref{eq:2.8} is optimal if and only if
\begin{equation}\label{eq:2.14}
M=A^{-1},\qquad
\alpha_r(\beta_r-g_r^Tc)-q_r^TAg_r=0\quad\text{for every }r.
\end{equation}

For each $\alpha_r>0$,
\begin{equation}\label{eq:2.15}
\beta_r-g_r^Tc=\|Ag_r\|,\qquad
q_r=\alpha_r\frac{Ag_r}{\|Ag_r\|}.
\end{equation}

With $p_r=Ag_r/\|Ag_r\|$ and $\mu_r=\alpha_r\|Ag_r\|$, the vectors and weights satisfy \eqref{eq:2.11}.

\end{proposition}

\begin{proof}
At a primal--dual optimum, Proposition~\ref{prop-exact-duality} gives
$\det A\det M=1$. Equality throughout the trace--AM--GM chain \eqref{eq:2.7} then forces
$Y=M^{1/2}AM^{1/2}=I$ and $\operatorname{tr}(AM)=n$.
The cone pairings are nonnegative and, by balance and normalization, their
sum is $n-\operatorname{tr}(AM)=0$. Thus every pairing vanishes.
Conversely, $M=A^{-1}$ gives equality in the certificate bound, proving
optimality of the feasible pair. For $\alpha_r>0$, equality in the support
and Cauchy--Schwarz inequalities gives \eqref{eq:2.15}.

The remaining point is to recover the full second moment from the symmetrized
aggregate. To recover isotropy, set $T=\sum_r\alpha_rg_rg_r^T/\|Ag_r\|$.
The aggregate equation is a Sylvester equation for $T$:
$AT+TA=2A^{-1}$. For $H=T-A^{-2}$, therefore,
\begin{equation}\label{eq:2.16}
AH+HA=0,
\qquad
0=\operatorname{tr}\bigl(H(AH+HA)\bigr)
=2\operatorname{tr}(AH^2)
\ge2\lambda_{\min}(A)\|H\|_F^2.
\end{equation}

Thus $H=0$ and $T=A^{-2}$. Positive definiteness has determined the
unsymmetrized moment matrix without a commutativity assumption.
The contact moments follow:
\begin{equation}\label{eq:2.17}
\sum_r\mu_rp_r=A\sum_r\alpha_rg_r=0,
\qquad
\sum_r\mu_rp_rp_r^T=ATA=I,
\qquad
\sum_r\mu_r=\operatorname{tr}I=n.
\end{equation}

This proves \eqref{eq:2.11}. 
\end{proof}

These results explain the exact optimization meaning of the certificate.
The universal contraction proof in Section~\ref{sec:3} uses Lemma~\ref{lem-certificate} and one explicit
feasible certificate; it does not require solving the normalized dual.
Section~\ref{sec:10.1} gives the usual unnormalized dual and its scale reduction.

\subsection{Iteration, rounding and computation}\label{sec:7}

The two geometric threads also clarify what the result contributes to a cutting-plane method. The contraction theorem gives an exact potential decrease. The certificate gives an instance-dependent bound and a stopping criterion for an ellipsoid subproblem. Rounding supplies a controlled coordinate system for subsequent computations.

\subsubsection{The exact-center iteration bound}\label{sec:7.1}

\begin{corollary}[{successive exact central cuts}]\label{cor-iteration}
Suppose
$K_{j+1}=K_j\cap H_j$, where each cutting hyperplane passes through $c(K_j)$. Then
\begin{equation}\label{eq:7.1}
w(K_k)\le R_n^k w(K_0).
\end{equation}

If $K_0$ is contained in a ball of radius $R$ and all retained regions contain a target ball of radius $r>0$, then
\begin{equation}\label{eq:7.2}
k\le\frac{n\log(R/r)}{-\log R_n}.
\end{equation}

\end{corollary}

\begin{proof}
Iterate the contraction, then compare the resulting upper bound
with the positive volume that the retained target ball requires.
Induction using Theorem~\ref{thm-sharp} gives \eqref{eq:7.1}. Monotonicity of $w$ and the two ball inclusions imply
\begin{equation}\label{eq:7.3}
\kappa_n r^n\le w(K_k)\le R_n^k w(K_0)\le R_n^k\kappa_n R^n.
\end{equation}

Since $\kappa_n>0$ and $0<R_n<1$,
\begin{equation}\label{eq:7.4}
n\log r\le k\log R_n+n\log R
\quad\Longrightarrow\quad
k(-\log R_n)\le n\log(R/r),
\end{equation}
which is \eqref{eq:7.2}. 
\end{proof}

The product bound treats successive cuts separately. A natural direction is to investigate whether compatibility of near-equality geometries across recentering yields a sharper bound for consecutive cuts.

Thus an oracle process satisfying the retained-ball assumption cannot continue producing such cuts beyond this bound. For minimization, an appropriate objective-sublevel inradius must be established before turning \eqref{eq:7.2} into an objective-accuracy guarantee. This qualification is the same one that enters other volume-based localization analyses.

The exact constant improves the finite-dimensional coefficient, while preserving the known $O(n\log(R/r))$ order of the inscribed-ellipsoid method. Its dimension dependence is explicit:
\begin{equation}\label{eq:7.5}
-\log R_n=\log2-\frac12+
\sum_{j=1}^{\infty}\frac1{2j(j+1)n^j}
=\log2-\frac12+\frac1{4n}+\frac1{12n^2}+\frac1{24n^3}+O(n^{-4}).
\end{equation}

For $n\ge2$, the absolutely convergent logarithmic series gives
\begin{equation}\label{eq:7.6}
\begin{aligned}
-\log R_n
&=\log2+\frac{n-1}{2}\log(1-1/n)\\
&=\log2-\frac12
+\sum_{j\ge1}\left(\frac1{2j}-\frac1{2(j+1)}\right)n^{-j}.
\end{aligned}
\end{equation}

This proves \eqref{eq:7.5}; its positive series terms decrease with $n$, so $R_n\uparrow\sqrt e/2$. Together with $R_1<R_2$, the monotonicity holds in every positive dimension.

\subsubsection{Rounding also helps initialize the next subproblem}\label{sec:7.2}

The inclusion \eqref{eq:1.2} makes the value of computing a John ellipsoid broader than center selection. In its coordinates, $B_2^n\subseteq K\subseteq nB_2^n$, so both a feasible inner ball and a controlled outer radius are available.

After a central cut, the ball with center $-v/4$ and radius $1/8$ remains feasible. Translating and scaling therefore give the explicit bound
\begin{equation}\label{eq:7.7}
B_2^n\subseteq 8\bigl((K\cap H_v)+v/4\bigr)
\subseteq(8n+2)B_2^n.
\end{equation}

For $\|z\|\le1$, the point $x=-v/4+z/8$ satisfies
\begin{equation}\label{eq:7.8}
\|x\|\le\frac38<1,\qquad v^Tx\le-\frac14+\frac18=-\frac18<0,
\end{equation}
proving the first inclusion. For $x\in K\cap H_v$, the rounding bound gives
$\|8(x+v/4)\|\le8\|x\|+2\le8n+2$, proving the second.
Thus an exact John normalization supplies an explicit next-subproblem initialization with radius ratio $O(n)$. For repeated computation, natural extensions include feasible warm starts and certified cuts through approximate centers.

\subsubsection{Computable bounds for individual instances}\label{sec:7.3}

The universal theorem bounds the worst case. A particular instance can admit a stronger certificate. One may fix the balanced weights in \eqref{eq:3.5} and optimize only the vector multipliers, or optimize all multipliers in the normalized dual. The explicit certificate is feasible for both choices, so such optimization can only improve its determinant bound.

For any polytope ellipsoid subproblem, a primal--dual feasible pair with gap $G$ satisfies
\begin{equation}\label{eq:7.9}
\operatorname{vol}\mathcal E(c,A)\ge e^{-G}w(\mathcal P).
\end{equation}

Indeed, with $p_*:=\log(w(\mathcal P)/\kappa_n)$, weak duality gives
\begin{equation}\label{eq:7.10}
p_*\le-\log\det M=\log\det A+G
\quad\Longrightarrow\quad
\frac{\operatorname{vol}\mathcal E(c,A)}{w(\mathcal P)}
=e^{\log\det A-p_*}\ge e^{-G}.
\end{equation}

This is a relative-volume stopping criterion for a primal--dual feasible computation. Relating this gap to center error could lead to certified contraction guarantees for cuts through computed centers.

For a fixed John frame and cut direction, let $\mathfrak D(\mathcal U,v)$ denote the maximum determinant in the normalized dual. Exact duality and the contact reduction give
\begin{equation}\label{eq:7.11}
\frac{w(P_{\mathcal U}\cap H_v)}{\kappa_n}
=\mathfrak D(\mathcal U,v)^{-1},\qquad
C_n=\left(\inf_{\mathcal U,v}\mathfrak D(\mathcal U,v)\right)^{-1},
\qquad \inf_{\mathcal U,v}\mathfrak D(\mathcal U,v)=R_n^{-1}.
\end{equation}

The last infimum is attained by the equality models. Exact duality and
contact recovery were proved in Section~\ref{sec:2.3}. Section~\ref{sec:10} discusses further
properties of the certificate and related matrix inequalities.

\subsection{Further dual and matrix properties}\label{sec:10}

The preceding sections establish the sharp contraction, equality, and stability
results. We now discuss extensions of the optimization and matrix viewpoints:
the unnormalized dual, the scope of the common rotation's optimality, an
alternative deficit identity, and the full family of coupled inequalities.
These results elaborate the earlier arguments; the core proofs are complete
without referring forward to this discussion.

\subsubsection{The unnormalized dual and its scale}\label{sec:10.1}

For the bounded full-dimensional polytope $\mathcal P$ of Section~\ref{sec:2}, retain
the primal formulation \eqref{eq:2.3}. The normalized certificate problem \eqref{eq:2.8} was
chosen for a direct trace bound. The usual Lagrangian dual explains this
normalization and is useful when working with unnormalized multipliers.

This standard convex formulation is discussed in \cite[Section 8.4]{BV2004} and in the maximum-volume ellipsoid literature \cite{ZG2003}. Self-duality of the Lorentz cone and the log-determinant conjugate give the dual
\begin{equation}\label{eq:10.1}
\begin{aligned}
\text{minimize}_{\alpha,q}\quad&\sum_r\alpha_r\beta_r-n-\log\det M\\
\text{subject to}\quad&\sum_r\alpha_rg_r=0,\quad
\|q_r\|\le\alpha_r,\quad
M=\sum_r\operatorname{sym}(q_rg_r^T)\succ0.
\end{aligned}
\end{equation}

Indeed, with $s=\sum_r\alpha_r\beta_r$ the Lagrangian is
\begin{equation}\label{eq:10.2}
\begin{aligned}
\mathcal L(A,c)&=\log\det A-\operatorname{tr}(AM)
-c^T\sum_r\alpha_rg_r+s,\\
\sup_{A\succ0,c}\mathcal L&=
\begin{cases}
s-n-\log\det M,&\sum_r\alpha_rg_r=0,\ M\succ0,\\
+\infty,&\text{otherwise}.
\end{cases}
\end{aligned}
\end{equation}

The finite supremum is attained at $A=M^{-1}$. This is the standard log-determinant conjugate calculation \cite[Chapters 3 and 5]{BV2004}.

\textbf{Remove the redundant multiplier scale.} Balance and cone feasibility
are preserved by positive scaling. We can therefore optimize that one scalar
before considering the remaining dual variables. For a dual feasible point,
primal feasibility implies $s\ge\operatorname{tr}(AM)>0$.
Scaling all multipliers by $\rho>0$ gives
$\Psi(\rho)=\rho s-n-n\log\rho-\log\det M$. With $u=\rho s/n$,
\begin{equation}\label{eq:10.3}
\Psi(\rho)-\Psi(n/s)=n(u-1-\log u)\ge0.
\end{equation}

Thus the optimal scale has $\sum_r\alpha_r\beta_r=n$, giving precisely
the normalized certificate problem \eqref{eq:2.8}. Its optimal value is the negative
of the primal optimal value, and both optima are attained, by Proposition~\ref{prop-exact-duality}.
A strictly feasible small ellipsoid also allows the use of Slater's theorem.
The construction in Section~\ref{sec:2.3} additionally identifies optimal multipliers
from the contact geometry.

\subsubsection{Optimality of the common rotation}\label{sec:10.2}

\begin{proposition}[{optimality among common contractions}]\label{prop-polar}
If $\|R\|_{\rm op}\le1$ and $N=\operatorname{sym}(RQ^T)\succ0$, then $\det N\le\det Q$. The rotation $R=U$ in \eqref{eq:3.13} attains equality.

\end{proposition}

\begin{proof}
We compare the determinant with its positive symmetric part,
then use the contraction bound on $R$. Write $RQ^T=N+K$ with $K^T=-K$,
and normalize the skew part by $T=N^{-1/2}KN^{-1/2}$.
The following Gram determinant identity provides the first comparison:
\begin{equation}\label{eq:10.4}
\det(I+T)^2=\det\bigl((I+T)^T(I+T)\bigr)
=\det(I+T^TT)\ge1
\end{equation}
gives $|\det(N+K)|\ge\det N$. Moreover,
\begin{equation}\label{eq:10.5}
\|R\|_{\rm op}\le1
\quad\Longrightarrow\quad 0\preceq R^TR\preceq I
\quad\Longrightarrow\quad |\det R|^2=\det(R^TR)\le1.
\end{equation}

Consequently $\det N\le|\det(RQ^T)|=|\det R|\det Q\le\det Q$.
At $R=U$, \eqref{eq:3.15} gives $N=M$ and equality. 
\end{proof}

This restricted optimality does not rule out improvements from changing the multiplier directions individually. The proof uses only real positive-semidefinite matrices, avoiding a separate argument with complex eigenvalues of a skew-symmetric matrix.

\subsubsection{An alternative deficit decomposition}\label{sec:10.3}

This decomposition follows the same gap through singular-value estimates.
Its terms distinguish support slack, the triangle inequality for the nuclear
norm, and the scalar logarithmic inequality. A fourth term records the
difference between the resulting nuclear-norm bound and the trace used by
the polar certificate.
Let $\|\cdot\|_*$ denote the sum of singular values, and put
\begin{equation}\label{eq:10.6}
\begin{aligned}
S&=\sum_i\alpha_i(1-c^Tu_i-\|Au_i\|)+t(-c^Tv-\|Av\|),\\
T_{\rm nuc}&=\sum_i\alpha_i\|Au_i\|+t\|Av\|-\|AQ\|_*,\\
H_{\rm nuc}&=\|AQ\|_*-n-\log\det(AQ),\\
\mathcal R_{\rm pol}&=\|AQ\|_*-\operatorname{tr}(AM).
\end{aligned}
\end{equation}

Primal feasibility gives $S\ge0$. For the other terms, use the standard rank-one and trace bounds for the nuclear norm \cite{HJ1985}. If $\rho_j>0$ are the singular values of $AQ$, then
\begin{equation}\label{eq:10.7}
\begin{aligned}
\|AQ\|_*
&\le\sum_i\alpha_i\|Au_i p_i^T\|_*+t\|Avv^T\|_*\\
&=\sum_i\alpha_i\|Au_i\|+t\|Av\|,\\
H_{\rm nuc}&=\sum_{j=1}^n(\rho_j-1-\log\rho_j)\ge0,\\
\operatorname{tr}(AM)&=\operatorname{tr}(AQU^T)
\le\|AQ\|_*\|U\|_{\rm op}=\|AQ\|_*.
\end{aligned}
\end{equation}

Here $\det(AQ)>0$ ensures $\sum_j\log\rho_j=\log\det(AQ)$.
All four terms are nonnegative. To relate them to the main deficit identity,
keep the polar correction explicitly: it is added to the cone gap and
subtracted from the matrix divergence, so it cancels in their sum.
Direct calculation gives
\begin{equation}\label{eq:10.8}
\begin{aligned}
\mathcal C&=S+T_{\rm nuc}+\mathcal R_{\rm pol},\\
\mathcal D(A,M)&=H_{\rm nuc}-\mathcal R_{\rm pol},\\
\epsilon&=\log\det\widehat F+S+T_{\rm nuc}+H_{\rm nuc}.
\end{aligned}
\end{equation}

In general $\mathcal C\ne S+T_{\rm nuc}$, because multiplication by $A$ changes the polar factor. The correction vanishes when $A$ commutes with $M$. The cone-based identity \eqref{eq:4.3} is more directly suited to equality and stability, because its two optimization terms separately measure complementarity and matrix matching.

\subsubsection{General and complementary coupled inequalities}\label{sec:10.4}

Lemma~\ref{lem-variational} extracted the particular coupling estimate used in the main proofs.
The same finite moment calculation gives an inequality for every symmetric
test matrix and every vector. We record that full statement and its companion
at the opposite endpoint of the trace range.

For these extensions introduce
\begin{equation}\label{eq:10.9}
\begin{gathered}
h=d+b^2,\qquad\sigma_*:=\frac{2(1-b)}h,\qquad\tau=\sigma_*-\sigma,\\
e_*:=\frac{2b}h=\frac{\sigma_*}{m},\qquad
\bar\chi=\sum_i\nu_i(1-s_i)(s_i+b)^2.
\end{gathered}
\end{equation}

Here $\tau$ is the unused part of the allowed trace, $e_*=\sigma_*/m$,
and $\bar\chi$ is the complementary nonnegative contact sum.
The coefficient in Lemma~\ref{lem-variational} is $L=d\tau/2$.
The scalar relation \eqref{eq:8.45} becomes
$\chi+(1-b)(2+b)^2\bar\chi/h^2=d(2+b)(1-b)/h$.

\begin{lemma}[{coupled Cauchy inequality}]\label{lem-cauchy}
For every symmetric $Z\in\mathbb S^m$ and $q\in\mathbb R^m$,
\begin{equation}\label{eq:10.10}
(\operatorname{tr}(ZE)+q^Tj)^2
\le\sigma\left[
\operatorname{tr}\bigl(Z^2(E+\tfrac12I_m)\bigr)
+2q^TZj+\frac d4\tau\|q\|^2\right].
\end{equation}

\end{lemma}

\begin{proof}
The aim is to represent the linear quantity on the left as one
weighted contact sum, so that one application of Cauchy--Schwarz controls both
parts. Choose $\gamma=\sqrt{d(1-b)}/h$; the moment identities then give
\begin{equation}\label{eq:10.11}
\operatorname{tr}(ZE)+q^Tj
=-\frac1d\sum_i\nu_i(s_i-b)y_i^T\bigl(Zy_i+\gamma(s_i+b)q\bigr).
\end{equation}

The coefficient $\gamma$ is chosen so that the $q$ term is exactly $q^Tj$.
Apply $|y_i^Tw|\le\|y_i\|\|w\|$ to each summand and use
$\|y_i\|^2=(1-s_i)(1+s_i)$. In weighted Cauchy--Schwarz, the factors
$\nu_i(1+s_i)(s_i-b)^2$ sum to $\chi$, leaving the bound
\begin{equation}\label{eq:10.12}
\frac{\chi}{d^2}\sum_i\nu_i(1-s_i)
\|Zy_i+\gamma(s_i+b)q\|^2.
\end{equation}

Terms at the poles have $y_i=0$, so no division by $1\pm s_i$ is needed.
Expanding the remaining squared norm produces a quadratic term in $Z$, a
mixed term in $Z,q$, and a quadratic term in $q$. The following three identities
evaluate them in that order:
\begin{equation}\label{eq:10.13}
\begin{aligned}
\sum_i\nu_i(1-s_i)y_iy_i^T&=\frac{2d}{2+b}(E+\tfrac12I_m),\\
\gamma\sum_i\nu_i(1-s_i)(s_i+b)y_i&=\frac{2d}{2+b}j,\\
\gamma^2\bar\chi&=\frac{d^2}{2(2+b)}\tau.
\end{aligned}
\end{equation}

The first two follow from \eqref{eq:8.44} and \eqref{eq:8.47}; the third uses the complementary
defect relation \eqref{eq:8.45}. Thus $\tau$ enters as the coefficient of the pure
$q$ term. Combining the three contributions gives
\begin{equation}\label{eq:10.14}
\begin{aligned}
&\sum_i\nu_i(1-s_i)\|Zy_i+\gamma(s_i+b)q\|^2\\
&\qquad=\frac{2d}{2+b}\left[
\operatorname{tr}\bigl(Z^2(E+\tfrac12I_m)\bigr)
+2q^TZj+\frac d4\tau\|q\|^2\right].
\end{aligned}
\end{equation}

Substitution and $\sigma=2\chi/[d(2+b)]$ prove \eqref{eq:10.10}. 
\end{proof}

\textbf{The complementary endpoint.} Replacing the transverse perturbation by
$e_*I_m-E$ replaces its trace $\sigma$ by $\tau$. The corresponding estimate is

\begin{lemma}[{complementary Cauchy inequality}]\label{lem-cauchy-minus}
With the notation of Lemma~\ref{lem-moments}, for every symmetric $Z$ and every $q$,
\begin{equation}\label{eq:10.15}
\bigl(\operatorname{tr}(Z(e_*I_m-E))+q^Tj\bigr)^2
\le\tau\left[\operatorname{tr}\bigl(Z^2(tI_m-E)\bigr)
+2q^TZj+\frac d4\sigma\|q\|^2\right].
\end{equation}

\end{lemma}

\begin{proof}
Use the complementary defect $\bar\chi$ in the same weighted
Cauchy--Schwarz argument as Lemma~\ref{lem-cauchy}. Choose
$c_-=(1-ab)/h$ and $\gamma_-=\sqrt d/[\sqrt{1-b}(2+b)]$ so that the
linear contact sum recovers both the shifted transverse block and the
coupling with the required coefficients. The moment identities give
\begin{equation}\label{eq:10.16}
\operatorname{tr}(Z(e_*I_m-E))+q^Tj
=c_-\sum_i\nu_i(s_i+b)y_i^T\bigl(Zy_i-\gamma_-(s_i-b)q\bigr).
\end{equation}

Weighted Cauchy--Schwarz bounds the square of \eqref{eq:10.16} by
\begin{equation}\label{eq:10.17}
c_-^2\bar\chi\sum_i\nu_i(1+s_i)
\|Zy_i-\gamma_-(s_i-b)q\|^2.
\end{equation}

The complementary defect relation turns the first factor into a multiple of
$\tau$. The expanded squared norm contains the opposite transverse moment
$tI_m-E$. Using \eqref{eq:8.47} and the moment identities gives
\begin{equation}\label{eq:10.18}
\begin{aligned}
&c_-^2\bar\chi=\frac{(1-b)(2+b)}{2d}\tau,\\
&\sum_i\nu_i(1+s_i)\|Zy_i-\gamma_-(s_i-b)q\|^2\\
&\qquad=\frac{2d}{(1-b)(2+b)}\left[
\operatorname{tr}\bigl(Z^2(tI_m-E)\bigr)
+2q^TZj+\frac d4\sigma\|q\|^2\right].
\end{aligned}
\end{equation}

The prefactors cancel to give \eqref{eq:10.15}. Contacts with $s_i=\pm1$ have $y_i=0$, so the argument covers the poles without division by zero. 
\end{proof}

Taking $Z=0$ in either \eqref{eq:10.10} or \eqref{eq:10.15} gives $\|j\|^2\le d\sigma\tau/4$.

\subsubsection{Inverse-trace interpretation}\label{sec:10.5}

\begin{lemma}[{inverse-trace interpretation}]\label{lem-inverse-trace}
Suppose $0<\sigma<\sigma_*$. Set $\eta=d\tau/4$ and $x=d\sigma/4$. Then \eqref{eq:10.10}, for all $Z,q$, is equivalent to
\begin{equation}\label{eq:10.19}
R_+=E+\frac12I_m-\frac{jj^T}{\eta}\succ0,
\qquad \operatorname{tr}(R_+^{-1})\le2m.
\end{equation}

Similarly, \eqref{eq:10.15} is equivalent to
\begin{equation}\label{eq:10.20}
R_-=tI_m-E-\frac{jj^T}{x}\succ0,
\qquad \operatorname{tr}(R_-^{-1})\le\frac{m}{t-e_*}.
\end{equation}

\end{lemma}

\begin{proof}
First absorb the mixed term by completing the square. Then identify
the best Cauchy--Schwarz coefficient and exhibit a test attaining it. This will
prove necessity as well as sufficiency of the inverse-trace bound.

\textbf{Step 1: separate the matrix and vector tests.} Complete the square in
\eqref{eq:10.10} by setting $q=q'-Zj/\eta$. This shift removes the mixed quadratic term
and leaves $R_+$ as the effective matrix coefficient. The two forms become
\begin{equation}\label{eq:10.21}
\begin{aligned}
\operatorname{tr}(ZE)+q^Tj
&=\operatorname{tr}\bigl(Z(R_+-\tfrac12I_m)\bigr)+q'^Tj,\\
\operatorname{tr}\bigl(Z^2(E+\tfrac12I_m)\bigr)+2q^TZj+\eta\|q\|^2
&=\operatorname{tr}(Z^2R_+)+\eta\|q'\|^2.
\end{aligned}
\end{equation}

Before using its inverse, establish strict positivity from the inequality
itself. If $w^TR_+w<0$, the choice $Z=ww^T$, $q'=0$ makes its right-hand side
negative. If $R_+w=0$ for $w\ne0$, the same choice gives zero on the right and
$(-\|w\|^2/2)^2>0$ on the left. Thus $R_+\succ0$.

\textbf{Step 2: identify the coefficient required for every test.}
Diagonalize $R_+$ and apply weighted Cauchy--Schwarz to the diagonal entries
of $Z$ and to $q'$; off-diagonal entries contribute only nonnegative terms to
$\operatorname{tr}(Z^2R_+)$. The resulting coefficient is
\begin{equation}\label{eq:10.22}
\begin{aligned}
S_+
&=\operatorname{tr}\bigl((R_+-\tfrac12I_m)^2R_+^{-1}\bigr)
+\frac{\|j\|^2}{\eta}\\
&=\operatorname{tr}R_+-m+\frac14\operatorname{tr}(R_+^{-1})
+\frac{\|j\|^2}{\eta}\\
&=\sigma-\frac m2+\frac14\operatorname{tr}(R_+^{-1}).
\end{aligned}
\end{equation}

To show that no smaller coefficient suffices, use the test attaining equality
in these scalar Cauchy--Schwarz inequalities:
\begin{equation}\label{eq:10.23}
Z_+=I_m-\frac12R_+^{-1},\qquad q'_+=j/\eta.
\end{equation}

At this test the linear form and the quadratic form both equal $S_+$.
The cancellation in \eqref{eq:10.22} leaves only the inverse trace to be bounded.
Consequently, the inequality for all $Z,q$ is equivalent to $S_+\le\sigma$, or
$\operatorname{tr}(R_+^{-1})\le2m$. This proves both directions of \eqref{eq:10.19}, including the case $S_+=0$.

\textbf{Step 3: apply the same argument at the complementary endpoint.}
For \eqref{eq:10.15}, set $q=q'-Zj/x$ and put
$\beta=t-e_*=(1-b)^2(2+b)/(2h)>0$. The completed forms have the same shape with
$(R_+,\tfrac12,\eta)$ replaced by $(R_-,\beta,x)$. The rank-one tests again give $R_-\succ0$. Weighted Cauchy--Schwarz is attained at
$Z_-=I_m-\beta R_-^{-1}$, $q'_-=j/x$, with coefficient
\begin{equation}\label{eq:10.24}
\begin{aligned}
S_-
&=\operatorname{tr}\bigl((R_--\beta I_m)^2R_-^{-1}\bigr)
+\frac{\|j\|^2}{x}\\
&=(-t+2e_*)m-\sigma+(t-e_*)^2\operatorname{tr}(R_-^{-1}).
\end{aligned}
\end{equation}

Since $\tau=me_*-\sigma$,
\begin{equation}\label{eq:10.25}
S_-\le\tau
\iff (t-e_*)^2\operatorname{tr}(R_-^{-1})\le m(t-e_*)
\iff \operatorname{tr}(R_-^{-1})\le\frac{m}{t-e_*}.
\end{equation}

This is \eqref{eq:10.20}, with both implications proved. 
\end{proof}

\subsubsection{Further models and extensions}\label{sec:10.6}

For a simplex with vertices $x_0,\ldots,x_n$ and centroid $\bar x$, the standard
affine-transport formula is $J(\Delta)=\bar x+H_\Delta^{1/2}B_2^n$, where
$H_\Delta=[n(n+1)]^{-1}\sum_i(x_i-\bar x)(x_i-\bar x)^T$.
It follows by transporting the regular-simplex inball. The main proof needs
only its contact configuration and the support calculation in Section~\ref{sec:3.2}.

For a full-dimensional offset-cut polytope, the same dual construction uses
the changed normalization below. For several cuts, the balance condition
contains all their normals:
\begin{equation}\label{eq:10.26}
\sum_i\alpha_i+\eta\alpha_0=n\quad(v^Tx\le\eta),\qquad
\sum_i\alpha_iu_i+\sum_\ell\alpha_\ell v_\ell=0.
\end{equation}

These formulations identify how the certificate changes under additional constraints. Quantitative bounds for offset cuts and successive recentered cuts are natural directions for further research, alongside exact offset profiles and sharp bounds for prescribed simultaneous cuts.

For completeness, the Lorentz transformation underlying \eqref{eq:3.7} is
\begin{equation}\label{eq:10.27}
\begin{pmatrix}1&-av^T\\-av&vv^T+kP\end{pmatrix}\binom1u
=\binom{1-as}{(s-a)v+kPu},\qquad s=u^Tv.
\end{equation}

Dividing the second component by the positive first component gives the unit
direction $p$ used in the certificate. This is the same construction in every
dimension; the dimensional distinction lies in its determinant proof.

\section{Conclusion}\label{sec:11}

The sharp central-cut factor is $R_n$ in every positive dimension, and equality is attained. An explicit certificate links the original John contacts to an upper bound on the new ellipsoid volume. The exact deficit identity then separates contact geometry, cone complementarity, and shape matching, yielding the complete equality family and quantitative stability.

The certificate construction is uniform in dimension, while its determinant proof still treats dimensions two, three, and at least four separately. A unified proof and sharper stability constants remain natural questions. Further directions include bounds for offset cuts, quantitative contraction over successive cuts, and certified cuts through approximate centers. Exact offset profiles, sharp simultaneous-cut bounds, and quantitative correction costs for computed centers remain open; whole-body stability requires hypotheses compatible with the nonunique equality family.

The results also clarify the role of inscribed ellipsoids in localization: contraction controls the potential, rounding controls a coordinate normalization, and a feasible primal--dual gap certifies the quality of an ellipsoid computation. These statements keep geometric progress separate from the cost and precision of centering. The analytical proofs are accompanied by a completed Lean 4 formalization.

\paragraph{Formal verification and reproducibility.}
The formalization uses Lean \texttt{v4.35.0-rc2} and mathlib at commit
\texttt{0a6c8e0355da0405d616f80b9f8232c4fab2cc5b} \cite{deMouraUllrich2021,Mathlib2020}.
It covers all 28 named results and their supporting mathematical statements,
with explicit correspondences for all 198 numbered equations. Independent
source builds and transitive axiom audits accept only \texttt{propext},
\texttt{Classical.choice}, and \texttt{Quot.sound}. The accompanying source
distribution, \texttt{SharpContractionCutJohnEllipsoid}, provides pinned
dependencies, statement maps, and reproducible verification scripts. The three
cited background assertions concerning centroid cuts and algorithmic complexity
in Section~\ref{sec:1.1} are outside the agreed formalization scope.

\section*{Use of AI-assisted tools}
The author used large language models to assist with the discovery of the
results, the development of the mathematical proofs, and the formal
verification in Lean 4.
\bibliographystyle{plainnat}
\bibliography{references}
\end{document}